\documentclass[
a4paper,reqno]{amsart}
\usepackage[T1,T5]{fontenc}
\usepackage[english]{babel}
\usepackage{amssymb,bbm,enumerate}
\usepackage{color}
\usepackage{mathtools}

\usepackage[
    linktocpage=true,
    colorlinks=true,
    linkcolor=blue,
    citecolor=red,
    urlcolor=green,
    pagebackref=true
]{hyperref}

\usepackage[all]{xy}
\usepackage{tikz-cd}

\usepackage{float}
\restylefloat{table}

\usepackage{tikz}

\usepackage{mathrsfs}
\usepackage[shortlabels]{enumitem}         
\usepackage{hyperref}         

\usepackage{ulem}

\newcommand{\id}{\mathrm{id}}

\DeclareMathOperator{\Hom}{Hom}
\DeclareMathOperator{\End}{End}

\DeclareMathOperator{\Ad}{Ad}

\DeclareMathOperator{\Der}{Der}

\newcommand{\Lie}{\mathcal{L}ie}

\newcommand{\g}{\mathfrak{g}}
\renewcommand{\O}{\mathcal{O}}
\newcommand{\eps}{\varepsilon}
\newcommand{\dd}{\mathrm{d}}
\newcommand{\tr}{\triangleright}
\newcommand{\TR}{\blacktriangleright}
\newcommand{\an}{_{\mathrm{an}}}

\theoremstyle{plain}
\newtheorem{theorem}{Theorem}[section]
\newtheorem{lemma}[theorem]{Lemma}
\newtheorem{proposition}[theorem]{Proposition}
\newtheorem{corollary}[theorem]{Corollary}

\theoremstyle{definition}
\newtheorem{definition}[theorem]{Definition}
\newtheorem{example}[theorem]{Example}

\theoremstyle{remark}
\newtheorem{remark}[theorem]{Remark}

\numberwithin{equation}{section}

\definecolor{light}{gray}{.9}

\begin{document}

\title{Affine post-group schemes and post-Lie algebras}

\author{Jiefeng Liu}
\address{School of Mathematics and Statistics, Northeast Normal University, Changchun, 130024, China}
\email{liujf534@nenu.edu.cn}

\author{Jiayao Ma}
\address{School of Mathematics and Statistics, Northeast Normal University, Changchun, 130024, China}
\email{jiayaoma@nenu.edu.cn}

\author{L\d\uhorn c Ta}
\address{Department of Mathematics, University of Pittsburgh, Pittsburgh, PA, USA}
\email{ldt37@pitt.edu}

\begin{abstract}
We develop a post-group analogue of the classical theory of affine group schemes, namely post-group objects in the category of affine schemes, which we call affine post-group schemes. Their coordinate rings form co-post-Hopf algebras, which admit a second Hopf algebra structure induced by the sub-adjacent affine group. 
Furthermore, we give functorial constructions of the post-Lie algebra of an affine post-group scheme and of the Rota--Baxter Lie algebra (of arbitrary weight) of an affine Rota--Baxter group scheme. 
Our constructions recover the analogous functors for algebraic post-Lie groups and Rota--Baxter Lie groups.
\end{abstract}

\maketitle

\section{Introduction}
A post-Lie algebra, a nonabelian generalization of a pre-Lie algebra, relates two Lie brackets on the same vector space via a bilinear operation. These algebraic structures naturally appear in numerical integration: classical B-series are based on pre-Lie algebras and the Butcher--Connes--Kreimer Hopf algebra, whereas Lie–Butcher series are based on post-Lie algebras and their universal enveloping algebras. The corresponding group-level structures are post-groups.

In a foundational work, Bai, Guo, Sheng and Tang \cite{BGST24} introduced the notions of post-groups, unifying several group structures appearing in diverse areas, including (Lie-)Butcher groups from numerical integration, \(\mathcal{P}\)-groups of operads, and braided groups for the Yang--Baxter equation (YBE). When equipped with a smooth structure, a post-group becomes a post-Lie group, whose differentiation yields a post-Lie algebra, providing a linear counterpart to these group structures. 

This paper addresses how to formulate the post-group structure in the category of affine schemes, called an affine post-group scheme, and how to differentiate it.

\subsection{Post-groups, post-Lie algebras and post-Hopf algebras}
The notion of a post-Lie algebra was introduced by Vallette in his study of Koszul duality of operads \cite{V07}.
Subsequently, Munthe-Kaas and Lundervold interpreted it differential-geometrically: a flat connection with constant torsion on a manifold \(M\) endows \(\mathfrak{X}(M)\) with the structure of a post-Lie algebra, related to numerical integration on manifolds \cite{MKL2013}. Moreover, the theory of regularity structures has been pioneered by Hairer \cite{Hairer14}. Recently, natural post-Lie algebraic structures on combinatorial objects such as decorated rooted trees and multi-indices were discovered by Bruned and Katsetsiadis \cite{BK23}. For more applications in this direction, see \cite{BST25, JZ25}.

The enveloping algebra of a post-Lie algebra always carries a \(D\)-algebra, which was introduced by Munthe-Kaas and Wright in \cite{MKW2008}. Recently, the notion of a post-Hopf algebra was introduced by Li, Sheng and Tang \cite{LST24}, with the universal enveloping algebra of a post-Lie algebra as the prototype. Subsequent work by Al-Kaabi, Ebrahimi-Fard and Manchon \cite{AEM24} further compared post-Hopf algebras with \(D\)-algebras, characterizing the former as ``good \(D\)-algebras'' since they are always equipped with a Grossman--Larson product. In \cite{BLS25}, the authors introduced post-Hopf algebroids and showed applications in geometric numerical integration on manifolds.
Post-Lie algebras have been studied from different points of view, for examples, Lax pairs \cite{BGN10}, factorization theorems \cite{EM18, EMM17}, as well as relations to post-(Lie) group \cite{AEM24, BGST24, GLSW26}. Applications of these post-structures also include Magnus expansions \cite{BGST24, EMQ25}, probability theory \cite{ER24, G25}, quiver-theoretical YBE \cite{STZ24}, affine feedback control systems \cite{EGV23} and multiple zeta values\cite{BK25}.

\subsection{The Yang--Baxter equation and Rota--Baxter operators}
It is well known that the YBE is a fundamental equation in mathematical physics, originating from the works of Yang \cite{Yang67} and Baxter \cite{Bax72}. The notion of Rota--Baxter operators can be traced back to G. Baxter's 1960 work on probability theory, where he considered such operators on associative algebras \cite{Baxter1960}. Rota--Baxter operators on Lie algebras are closely related to solutions of the (modified) classical YBE \cite{S1983}, with the former serving as the operator form of the latter. 
Bai, Guo and Ni showed that Rota-Baxter Lie algebras naturally gave rise to post-Lie algebras\cite{BGN10}.

Solutions to the YBE have been extensively studied, and a major step forward was Drinfeld's proposal to investigate set-theoretical solutions \cite{Dri92}. Further studies were carried out in \cite{ESS99, LYZ00}. An important algebraic structure arising from this study is the notion of (skew) braces for constructing set-theoretical solutions to the YBE \cite{CJO14, GV17, Rump06}. Rota--Baxter (Lie) groups, introduced by Guo, Lang and Sheng \cite{GLS21}, differentiate at the identity to Rota–Baxter Lie algebras, and their relations with skew braces were first established in \cite{BG22}.
The study of solutions to the YBE in the setting of affine schemes goes back to Angiono, Galindo and Vendramin \cite{AGV17}, who introduced affine (skew) braces together with their coordinate rings, called Hopf co-braces. They did not directly investigate the YBE in this setting, but they did mention this idea.
This idea was subsequently pursued in \cite{MZL26}, where affine Rota–Baxter groups and affine skew braces were developed, while earlier and independently the third author of the present paper \cite{Ta25} named such solutions affine Yang–Baxter schemes.

\subsection{Motivations and main results}

The (formal) integration problem for (complete) post-Lie algebras and (complete) Rota--Baxter Lie algebras has been studied successively \cite{BGST24, GKST26, GLS21, MQS20, MQ21}. We focus on the Lie group level, where differentiating a post-Lie group or a Rota--Baxter Lie group at the identity yields, respectively, a post-Lie algebra or a Rota--Baxter Lie algebra; see \cite[Section 4]{BGST24} and \cite[Section 2.2]{GLS21}. Clearly, these constructions are analytic and hence restricted to $k=\mathbb{R}$ or $\mathbb{C}$; the corresponding Lie algebras are also finite-dimensional.

To relax these restrictions on the base field $k$ and on finite-dimensionality, we work with group schemes: we define affine post-group schemes, study their differentiation, and construct a functor from the category of affine post-group schemes to that of post-Lie algebras (Theorems \ref{Thm the post-Lie alg of affine post-group} and \ref{Thm functor affine post group to postlie}); we also differentiate affine Rota--Baxter group schemes and construct a functor from their category to that of Rota–Baxter Lie algebras (Theorems \ref{Thm Rota--Baxter Lie algebra of affine RB group} and \ref{Thm functor affine RB group to RBlie}). In Section \ref{sec:comp}, we compare our algebraic constructions with analytic constructions (Theorems \ref{thm:post-lie} and \ref{thm:rb-lie}).

On the one hand, although affine skew braces were introduced in \cite{AGV17}, the corresponding post-group structure in the affine scheme setting does not appear to have been investigated. On the other hand, the dual notion of a post-Hopf algebra remains absent from the literature. Although Al-Kaabi, Ebrahimi-Fard and Manchon mentioned in \cite[Remark 7]{AEM24} that such a dual notion could be obtained by dualizing the axioms, to the best of our knowledge, no definition has been given and no further study has been carried out.
These two gaps partly motivate the present paper. We show that they are in fact two sides of the same coin (Theorems \ref{Thm affine post-group equivalent to co-post Hopf} and \ref{Thm AffPGrp simeq CoHopf opp}): the coordinate ring of an affine post-group scheme is precisely the dual notion to that of a post-Hopf algebra, which we call a co-post-Hopf algebra. The main property of a co-post-Hopf algebra is that one can associate to it a second Hopf algebra, called its sub-adjacent Hopf algebra (Proposition \ref{Prop sub adjacent Hopf algebra} and Remark \ref{Remark sub=sub}). The name ``sub-adjacent Hopf algebra'' reflects the fact that this Hopf algebra arises from the sub-adjacent affine group of an affine post-group.

In the context of the YBE, post-groups are known to provide set-theoretic solutions. We extend the set-theoretic construction from post-groups to affine schemes by showing that affine post-group schemes give rise to affine Yang–Baxter schemes. 

Another motivation comes from the categorical perspective. Group objects can be defined in any category with finite products, and this notion has proved remarkably fruitful: affine (super)group schemes and Lie (super)groups all arise as group objects in suitable categories.
This suggests that the notion of post-groups should also admit a categorical generalization. Accordingly, we introduce the notion of a post-group object, which simultaneously recovers post-groups (in the category of sets), post-Lie groups (in the category of smooth manifolds), and affine post-group schemes (in the category of affine schemes). 
We hope that this framework will extend to super analogues and other variants, and inspire further work.

The paper is organized as follows. In Section~\ref{sec pre}, we recall the necessary background on post-groups, affine group schemes, and introduce post-group objects as a categorical generalization of post-groups. In Section~\ref{sec affinepostgroup copostHopf}, we introduce affine (algebraic) post-group, prove that they admit a sub-adjacent affine group structure, and define co-post-Hopf algebras as their coordinate rings; as a consequence, a co-post-Hopf algebra naturally carries a second Hopf algebra structure. In Section~\ref{sec affine post-group schemes and related affine schemes}, we prove that the categories of affine post-groups and affine skew braces are isomorphic, which naturally yields affine Yang–Baxter schemes. In Section~\ref{sec post-Lie and RB Lie}, we construct functors that assign to each affine post-group scheme a post-Lie algebra, and to each affine Rota--Baxter group scheme a Rota--Baxter Lie algebra. 
In Section~\ref{sec:comp}, we compare our algebro-geometric constructions with the differential-geometric constructions in \cite{BGST24} and \cite{GLS21}; see Theorems \ref{thm:post-lie} and \ref{thm:rb-lie}.

Throughout, \(k\) is a field unless otherwise stated. All $k$-algebras are associative and unital and all algebra homomorphisms preserve the unit. In this paper, the terms co-post-Hopf algebra, co-Rota--Baxter Hopf algebra, and Hopf co-brace all refer to the commutative case, i.e., their underlying Hopf algebras are commutative.

\section{Preliminaries} \label{sec pre}
\subsection{Notation}
Let \(\mathsf{Set}\) denote the category of sets and \(\mathsf{Alg}_k\) the category of commutative \(k\)-algebras. For any category \(\mathcal{C}\), we write \(\mathcal{C}^{\mathrm{opp}}\) for its opposite category, \(\mathsf{Func}[\mathcal{C}, \mathsf{Set}]\) for the category of functors from \(\mathcal{C}\) to \(\mathsf{Set}\), and \(\mathsf{RepFunc}[\mathcal{C}, \mathsf{Set}]\) for its full subcategory of representable functors.

For convenience, we shall always assume that \(\mathcal{C}\) is a locally small category with finite products and coproducts and terminal object $*$.
Given two objects $A_1,A_2$, let \(\pi_i\colon A_1 \times A_2 \to A_i\)  
and \(\iota_i\colon A_i \to A_1 \otimes A_2\) denote the canonical morphisms.  
Given morphisms \(f_i\colon X \to A_i\), consider the induced morphism \(\langle f_1, f_2 \rangle\colon X \to A_1 \times A_2\). 
Given morphisms $h_i\colon A_i\to B_i$ with $i=1,2$, denote $h_1\times h_2\coloneq\langle h_1 \circ \pi_1, \, h_2 \circ \pi_2 \rangle$.
Write $A^{\times n}$ for the $n$-fold product. 
We denote $h^A\coloneq\Hom_\mathcal C(A,-)$, usually with $\mathcal C=\mathsf{Alg}_k$.

\subsection{Post-groups}
We first recall the definition of a post-group and introduce its basic properties.

\begin{definition}[(\cite{BGST24})]
\label{Def:postgroup}
A \textbf{post-group} is a group \((G, \cdot)\) equipped with a multiplication \(\rhd\) on \(G\) such that
\begin{enumerate}[label=(PG\arabic*), ref=(PG\arabic*)]
    \item \label{item:postgroup1} for each \(a \in G\), the left multiplication \(L_a^\rhd : G \to G\), \(L_a^\rhd(b) = a \rhd b\) for all \(b \in G\), is an endomorphism of the group \((G, \cdot)\), that is,
    \begin{equation*}
    a \rhd (b \cdot c) = (a \rhd b) \cdot (a \rhd c) \quad \text{for all } a, b, c \in G;
    \end{equation*}

    \item \label{item:postgroup2} for each $a \in G$, the map $L_a^{\rhd}$ is bijective;

    \item \label{item:postgroup3}the following “weighted” associativity holds:
    \begin{equation*}
    \big(a\cdot(a\rhd b)\big)\rhd c =a \rhd (b \rhd c)  \quad \text{for all } a, b, c \in G.   
    \end{equation*}
\end{enumerate}
Moreover, a post-group $(G, \cdot, \rhd)$ is called a \textbf{pre-group} if the group $(G, \cdot)$ is an abelian group.
\end{definition}

A \textbf{homomorphism} of post-groups is a group homomorphism preserving $\rhd$. Their category is denoted by $\mathsf{PGrp}$. The first two axioms say that each $L^\rhd_a$ is a group automorphism. The third explains why the following second group law appears naturally.

\begin{proposition}{\rm (\cite{BGST24})}
\label{prop post-group}
    Let $(G, \cdot, \rhd)$ be a post-group with identity $e$. Then $a \rhd e = e$ and $e \rhd a = a$. The operations 
    \begin{equation*}
      a \circ b=a \cdot (a\rhd b), \quad a^{\circ-1}=(L^{\rhd}_a)^{-1} (a^{-1})
    \end{equation*}
define a group $(G,\circ)$, which is called the \textbf{sub-adjacent group}, with the same identity. Moreover, $L^{\rhd}_{a \circ b}= L^{\rhd}_a L^{\rhd}_b$.
\end{proposition}

\begin{remark}
\label{Re a-circ-inv}
The inverse of \(a\) with respect to \(\circ\) can also be expressed as
\[
a^{\circ(-1)} = (L^{\rhd}_a)^{-1}(a^{-1}) = L^{\rhd}_{a^{\circ(-1)}}(a^{-1}).
\]
Note that $(G,\circ)$ is also called the \textbf{Grossman–Larson group} in \cite{AEM24, EGV23}.
\end{remark}

The following examples correspond to the constructions in
\cite[Section~5.1]{AEM24}.
\begin{example}
\label{ex:basic-post-groups}
Every group $(G,\cdot)$ admits the following post-group structures:
\begin{enumerate}
    \item[(i)] $a\rhd b=b$, with $a\circ b=a\cdot b$. Its sub-adjacent group coincides with $(G,\cdot)$.

    \item[(ii)] $a\rhd b=a^{-1}\cdot b\cdot a$, with $a\circ b=b\cdot a$. Its sub-adjacent group is $G^{\mathrm{op}}$, the same underlying set with reversed multiplication.
\end{enumerate}
\end{example}

\subsection{From post-groups to post-group objects}

\begin{definition}\label{Def postgroup object}
    A \textbf{post-group object in \(\mathcal{C} \)} is a quintuple \((G, m, e, \operatorname{inv}, \rhd)\), where \((G, m, e, \operatorname{inv})\) is a group object in \(\mathcal{C} \) and \(\rhd: G \times G \to G\) is a morphism satisfying the following conditions:
    \begin{enumerate}[label=(PGO\arabic*), ref=(PGO\arabic*)]
        \item \textbf{(Distributivity).}\label{item postgroup object 1}
        The following diagram commutes:
        \begin{equation}\label{Diag: postgroup object preserve multiplication}
        \begin{tikzcd}[column sep = huge, row sep = huge]
            G^{\times 3} \ar[r, "\operatorname{id} \times m"] \ar[d, "{\langle \operatorname{id}, \operatorname{id} \rangle \times \operatorname{id} \times \operatorname{id}}"] 
            & G^{\times 2} \ar[r, "\rhd"] & G \\
            G^{\times 4} \ar[r, "\operatorname{id} \times \tau \times \operatorname{id}"] 
            & G^{\times 4} \ar[r, "\rhd \times \rhd"] & G^{\times 2} \ar[u, "m"] ,
        \end{tikzcd}
        \end{equation}
        where \(\tau: G \times G \to G \times G\) is the swap, i.e., \(\tau = \langle\pi_2, \pi_1 \rangle\).

        \item \textbf{(Bijectivity).} \label{item postgroup object 2}
        For every object \(X \in \mathcal{C} \) and every morphism \(f\colon X \to G\), the induced map
        \[
        L^{\rhd}_{f}: \operatorname{Hom}_{\mathcal{C} }(X, G) \to \operatorname{Hom}_{\mathcal{C} }(X, G),\qquad
        g \mapsto \rhd \circ \langle f, g \rangle
        \]
        is a bijection.

        \item \textbf{(Weighted associativity).}
        \label{item postgroup object 3}The following diagram commutes:
        \begin{equation}
        \label{Diag: postgroup object weighted assoc}
        \begin{tikzcd}[column sep = huge, row sep = huge]
            G^{\times 3} \ar[r, "\operatorname{id} \times \rhd"] \ar[d, "{\langle \operatorname{id}, \operatorname{id} \rangle \times \operatorname{id} \times \operatorname{id}}"] 
            & G^{\times 2} \ar[r, "\rhd"] & G \\
            G^{\times 4} \ar[r, "\operatorname{id} \times \rhd \times \operatorname{id}"] 
            & G^{\times 3} \ar[r, "m \times \operatorname{id}"] & G^{\times 2} \ar[u, "\rhd"] .
        \end{tikzcd}
        \end{equation}
    \end{enumerate}
    For brevity, we often denote a post-group object by the triple \((G, m, \rhd)\). A post-group object $(G,m,\rhd)$ is called a {\bf pre-group object} when $(G,m)$ is an abelian group object.
\end{definition}

The notion of a post-group object in a category is a natural and reasonable generalization of the notion of a post-group.


\begin{proposition}
\label{Prop:postgroup object in Set= postgroup}
    When \(\mathcal{C}  = \mathsf{Set}\), the notion introduced above recovers exactly the notion of post-groups. That is, a post-group object in \(\mathsf{Set}\) is exactly a post-group.
\end{proposition}

\begin{proof}
    The isomorphism $\Hom(*,G)\cong G$ is compatible with $L^\tr$, so taking $X=*$ in \ref{item postgroup object 2} shows that each $L^\tr_a$ is bijective. 
    The rest of the proof is routine.
\end{proof}

\begin{example}
As defined in \cite[Definition 4.2]{BGST24}, a post-Lie group is a post-group \((G, \cdot, \triangleright)\) where \((G, \cdot)\) is a Lie group and \(\triangleright\) is smooth. Equivalently, a post-group object in the category of smooth manifolds is exactly a post-Lie group.
\end{example}

\subsection{Affine (algebraic) group schemes}
\label{sec affine algebraic group scheme}
We use the functor-of-points approach to discuss affine group schemes. For a detailed account, we refer the reader to Milne's book \cite{Milne17} and Waterhouse's book \cite{W1979}.

An \textbf{affine scheme} is a representable functor $G\colon \mathsf{Alg}_k \to \mathsf{Set}$. If $G$ factors through the forgetful functor $\mathsf{Grp}\to\mathsf{Set}$, we call $G$ an \textbf{affine group scheme} or simply an \textbf{affine group}. 
Equivalently, an affine group is a group object $(G,m,e, \operatorname{inv})$ in the category $\mathsf{RepFunc}[\mathsf{Alg}_k, \mathsf{Set}]$. 
We say that $G$ is \textbf{algebraic} or of \textbf{finite type} if its coordinate ring $\mathcal{O}(G)$ is finitely generated.

The usual equivalence of categories $\O\colon\mathsf{Aff}_k^{\mathrm{opp}}\xrightarrow{\sim}\mathsf{Alg}_k$ induces a contravariant equivalence of categories from affine group schemes to commutative Hopf algebras $(H,\mu,\eta,\Delta,\eps,S)$.
In calculations involving the comultiplication $\Delta$, we use Sweedler notation, suppressing summation. For example,
\[
\Delta(h)=h_1 \otimes h_2, \quad (\Delta \otimes \operatorname{id}) \circ \Delta(h)=h_1 \otimes h_2 \otimes h_3.
\]

\section{Affine post-group schemes and co-post-Hopf algebras} \label{sec affinepostgroup copostHopf} In this section, we introduce affine post-groups and co-post-Hopf algebras---the latter being the coordinate rings of the former. We prove that every affine post-group admits a new affine group structure, its sub-adjacent affine group, which in turn implies that a co-post-Hopf algebra carries a second Hopf algebra structure.

\subsection{Affine post-groups and sub-adjacent affine groups}
\begin{definition}
\label{Def affine post-group}
    An \textbf{affine post-group scheme over $k$} (or simply an \textbf{affine post-group}) is a representable functor $G\colon \mathsf{Alg}_k \to \mathsf{Set}$ with two natural transformations $m , \rhd \colon G\times G \to G$ such that, for any $R \in \mathsf{Alg}_k$,
    \[
        m_R\colon G(R)\times G(R) \to G(R) \quad \text{and} \quad \rhd_R\colon G(R) \times G(R) \to G(R)
    \]
   make $\big(G(R),m_R,\rhd_R \big)$ into a post-group, where $m_R$ is the group multiplication. If $G$ is represented by a finitely generated $k$-algebra, then it is called an \textbf{affine algebraic post-group scheme over k}.
\end{definition}

Equivalently, an affine post-group $(G,m,\rhd)$ is a functor $G:\mathsf{Alg}_k\to\mathsf{PGrp}$ whose underlying set-valued functor is representable. 


\begin{definition}
Let \((G, m_G, \rhd_G)\) and \((H, m_H, \rhd_H)\) be affine post-groups. A \textbf{homomorphism} \(\Phi\colon G \to H\) is a natural transformation such that, for any \(R \in \mathsf{Alg}_k\), the map $\Phi_R\colon G(R) \to H(R)$ is a homomorphism of post-groups. Equivalently, the following diagrams commute:
\begin{equation*}
\begin{tikzcd}
G \times G \arrow[r, "m_G"] \arrow[d, "\Phi \times \Phi"'] & G \arrow[d, "\Phi"] \\
H \times H \arrow[r, "m_H"'] & H
\end{tikzcd}
\qquad
\begin{tikzcd}
G \times G \arrow[r, "\rhd_G"] \arrow[d, "\Phi \times \Phi"'] & G \arrow[d, "\Phi"] \\
H \times H \arrow[r, "\rhd_H"'] & H .
\end{tikzcd}
\end{equation*}
\end{definition}

Affine post-groups over $k$ and their homomorphisms form a category, denoted by $\mathsf{AffPGrp}_k$. 

\begin{remark}
Let $(G,m,\rhd)$ be an affine post-group, and let $e$ and $\operatorname{inv}$ be the identity and inverse morphisms of $(G,m)$. For any  $R \in \mathsf{Alg}_k$ and $x,y \in G(R)$, denote
$m_R(x,y)=x \cdot_R y$, $\operatorname{inv}_R(x)=x^{-1}$ and $\rhd_R(x,y)=x \rhd_R y$. The identity element $1_R$ of $G(R)$ is given by the unique image of the map $(e \circ!)_R\colon G(R) \to G(R)$, where $!$ is the unique morphism from $G$ to the terminal object $h^k= \operatorname{Hom}_{\mathsf{Alg}_k}(k, -)$.
Namely, 
\[
1_R= e_R(\eta_R)=\eta_R \circ \varepsilon_{\mathcal{O}(G)}.
\] 
Here $\eta_R$ is the unit of $R$ and $\varepsilon_{\mathcal{O}(G)}$ is the counit of $\mathcal{O}(G)$.
\end{remark}

\begin{example}
\label{ex:affine-post-basic}
Let $(G,m)$ be an affine group scheme. The constructions in Example~\ref{ex:basic-post-groups} extend to affine post-groups as follows.
Let $R\in\mathsf{Alg}_k$.
\begin{enumerate}
    \item[(i)]
    The second projection 
    \[
    \rhd\coloneq \pi_2:G\times G\to G
    \]
    defines an affine post-group $(G,m,\rhd)$. Indeed, $x\rhd_R y=y$ for all $x,y\in G(R)$, so each left multiplication is the identity and the post-group identities hold immediately.

    \item[(ii)]
    The maps
    \[
     \rhd_R\colon G(R)\times G(R)\to G(R),\qquad (f,g) \mapsto f^{-1}\cdot_R g\cdot_R f,
    \]
    are natural in $R$, being defined by multiplication and inversion. They define an affine post-group $(G,m,\rhd)$ by Example~\ref{ex:basic-post-groups}\textup{(ii)}.
\end{enumerate}
\end{example}

See Example \ref{ex non-reduced example} for an example of an affine post-group whose underlying affine group is non-reduced.

The coordinate ring $\O(G)$ of an affine post-group $G$ is, first of all, a commutative Hopf algebra. Moreover, we will show in the next subsection that it carries the additional structure of a commutative {\bf co-post-Hopf algebra}; see Subsection \ref{subsection coordinate ring of affine post-group}.

\begin{lemma}
\label{Lem two equivalent bijective condition}
Let $(G,m,e,\operatorname{inv},\rhd)$ be an affine post-group. Then the following statements are equivalent.
\begin{enumerate}
    \renewcommand{\labelenumi}{\rm(\arabic{enumi})}
    \item \label{item bijection equivalence 2} For any $R \in \mathsf{Alg}_k$ and $a \in G(R)$, the map
    \[
    L_a^{\rhd_R}: G(R) \to G(R), \quad b \mapsto a\rhd_R b
    \]
    is bijective.
    \item \label{item bijection equivalence 1} For any $X \in \mathsf{RepFunc}[\mathsf{Alg}_k, \mathsf{Set}]$ and $\alpha \in \operatorname{Hom}(X,G)$, the map
    \[
    L_\alpha^\rhd: \operatorname{Hom}(X,G) \to \operatorname{Hom}(X,G), \quad \beta \mapsto \rhd \circ \langle \alpha, \beta \rangle
    \]
    is bijective.
\end{enumerate}
\end{lemma}

\begin{proof}
Since $G$ is an affine post-group, (1) holds by definition. We show that (1) and (2) are equivalent. We may assume that $X= h^R$ with $R\in\mathsf{Alg}_k$. Consider the Yoneda bijection
\[
\Phi_R\colon  \operatorname{Hom}(h^R, G) \xrightarrow{\cong} G(R), \qquad \Phi_R(\alpha) = \alpha_R(\operatorname{id}_R).
\]
For all $\alpha, \beta: h^R \to G$, we have
\[
\Phi_R\big(\rhd \circ \langle \alpha, \beta \rangle \big) = \big(\rhd \circ \langle \alpha, \beta\rangle\big)_R(\operatorname{id}_R) = \rhd_R\big(\alpha_R(\operatorname{id}_R), \beta_R(\operatorname{id}_R)\big) = \Phi_R(\alpha) \rhd_R \Phi_R(\beta),
\]
so 
\begin{equation}\label{diag L}
    \Phi_R \circ L^{\rhd}_{\alpha}=L_{\Phi_R(\alpha)}^{\rhd_R} \circ \Phi_R.
\end{equation}
Since $\Phi_R$ is bijective, $L_\alpha^\rhd$ is bijective if and only if $L_{\Phi_R(\alpha)}^{\rhd_R}$ is bijective. Since $\Phi_R$ is surjective, the proof is complete.
\end{proof}

\begin{proposition}\label{Prop: affine postgroup=postgroup object in repfunc}
    A quintuple \((G, m, e, \operatorname{inv}, \rhd)\) is an affine post-group over \(k\) if and only if 
    it is a post-group object in \(\mathsf{RepFunc}[\mathsf{Alg}_k, \mathsf{Set}]\).
\end{proposition}
\begin{proof}
For an affine post-group \((G, m, e, \operatorname{inv}, \rhd)\), the pair $(G,m)$ is an affine group, which is equivalent to a group object in \(\mathsf{RepFunc}[\mathsf{Alg}_k, \mathsf{Set}]\). 
Since equalities of natural transformations $\mathsf{Alg}_k\to\mathsf{Set}$ can be checked on $R$-points for each $R\in\mathsf{Alg}_k$, \ref{item:postgroup1} and \ref{item:postgroup3} hold for all $R\in\mathsf{Alg}_k$ if and only if \ref{item postgroup object 1} and \ref{item postgroup object 3} hold, respectively. 
The same goes for \ref{item:postgroup2} and \ref{item postgroup object 2} due to Lemma \ref{Lem two equivalent bijective condition}.
\end{proof}

\begin{proposition}
\label{prop:post-group}
Let \((G, m, e, \operatorname{inv}, \rhd)\) be an affine post-group, and let \(!: G \to h^k\) be the unique morphism to the terminal object. 
Then the following diagrams commute:
\begin{equation}
\label{Diag f rhd 1=1 and 1 rhd f = f}
\begin{tikzcd}[column sep =1.5cm]
G \ar[r, "e \circ !"] \ar[d, "{\langle \operatorname{id}, e\circ ! \rangle}"] & G \ar[equal,d]\\
G^{\times 2} \ar[r,"\rhd"]&  G 
\end{tikzcd}
\quad \text{and} \quad
\begin{tikzcd}[column sep =1.5cm]
G \ar[equal, r] \ar[d,"{\langle !,\operatorname{id} \rangle}"]& G\\
* \times G \ar[r,"e \times \operatorname{id}"]& G^{\times 2} \ar[u,"\rhd"].
\end{tikzcd}
\end{equation}
\end{proposition}
\begin{proof}
Fix $R\in\mathsf{Alg}_k$. According to Proposition \ref{prop post-group}, we have
\[
x \rhd_R 1_R=1_R \quad \text{and} \quad 1_R \rhd_R x=x,
\]
for any $x \in G(R)$,
which imply that 
\[
\big(\rhd\circ \langle \operatorname{id},e \circ !\rangle \big)_R= (e\circ!)_R \quad \text{and} \quad \big(\rhd \circ \langle e\circ!,\operatorname{id}\rangle \big)_R=\operatorname{id}_R.
\]
This completes the proof.
\end{proof}

An affine post-group naturally gives rise to a new affine group.

\begin{proposition}
\label{Prop sub-adjacent affine group of affine postgroup}
Let \((G, m, e, \operatorname{inv}, \rhd)\) be an affine post-group over \(k\). Define a natural transformation \(m^\rhd: G \times G \to G\) by the following commutative diagram:
\begin{equation*}
    \begin{tikzcd}[column sep = huge]
        G^{\times 2} \ar[r, "{m^{\rhd}}"] \ar[d, "{\langle \operatorname{id}, \operatorname{id} \rangle \times \operatorname{id}}"] & G 
        \\
        G^{\times 3} \ar[r, "\operatorname{id} \times \rhd"] & G^{\times 2} \ar[u, "m"].
    \end{tikzcd}
\end{equation*}
That is, for every $R\in \mathsf{Alg}_k$, 
\[
m_R^\rhd(x,y)=x\cdot_R(x\rhd_R y),
\]
and we write \(x\circ_R y\coloneq m_R^\rhd(x,y)\). Then \((G,m^\rhd,e,\operatorname{inv}^\rhd)\) is an affine group, where \(\operatorname{inv}^\rhd\colon G\to G\) is defined pointwise by
\[
  \operatorname{inv}^\rhd_R(x) = (L_x^{\rhd_R})^{-1}(x^{-1}).
\]
We denote this affine group by \(G_{\rhd}\) and call it the \textbf{sub-adjacent affine group} of the affine post-group \((G,m,e,\operatorname{inv},\rhd)\).
\end{proposition}
\begin{proof}
First, by Proposition~\ref{prop post-group}, the operation
\[
m^{\rhd}_R(x,y)=x \circ_R y = x \cdot_R (x \rhd_R y)
\]
defines a sub-adjacent group of the post-group $\big(G(R), m_R, \rhd_R \big)$ with the unit $1_R=e_R(\eta_R)$, and the inverse of $x$ in the sub-adjacent group is
\[
x^{\circ(-1)}=(L_x^{\rhd_R})^{-1}(x^{-1}).
\]
Thus, for each $R \in \mathsf{Alg}_k$, the functor $G$ makes $G(R)$ into a group $\big(G(R), \circ_R\big)$.
The naturality of \(m^\rhd\) and \(e\) is clear; it remains to check the naturality of \(\operatorname{inv}^\rhd\). Recall that $G(\varphi)$ is homomorphism of post-groups for any algebra morphism \(\varphi\colon  R \to S\). 
Moreover, by \cite[Theorem 2.4]{BGST24}, $G(\varphi)$ is also a homomorphism of sub-adjacent groups, which implies $G(\varphi)(x^{\circ(-1)})=G(\varphi)(x)^{\circ(-1)}$. 
Namely, 
\[
G(\varphi)\big(\operatorname{inv}_R^{\rhd}(x)\big)=\operatorname{inv}_S^{\rhd}\big(G(\varphi)(x) \big),
\]
which implies the naturality of $\operatorname{inv}^{\rhd}$.
Thus, \((G, m^\rhd, e, \operatorname{inv}^\rhd)\) is an affine group.
\end{proof}

\begin{example}
For the structures in Example~\ref{ex:affine-post-basic}:
\begin{enumerate}
    \item[(i)] $x\circ_R y=x\cdot_R y$, so the sub-adjacent affine group coincides with $(G,m)$.

    \item[(ii)] $x\circ_R y=y\cdot_R x$, so the sub-adjacent affine group is $G^{\mathrm{op}}\coloneq(G,m\circ\tau)$, where $\tau$ interchanges the factors.
\end{enumerate}
\end{example}

\begin{remark}
\label{Remark on sub-adjacent affine group}
The Proposition \ref{Prop sub-adjacent affine group of affine postgroup} implies that, for any $R \in \mathsf{Alg}_k$ and $x \in G(R)$,
    \[
    \begin{aligned}
        \operatorname{inv}_R^{\rhd}(x)=x^{\circ(-1)}&=(L_x^{\rhd_R})^{-1}(x^{-1})\\
        &=L_{x^{\circ(-1)}}^{\rhd_R}(x^{-1})\\
        &=x^{\circ(-1)} \rhd_R x^{-1}\\
        &=\rhd_R \circ \langle \operatorname{inv}_R^{\rhd},\operatorname{inv}_R \rangle (x) = \big(\rhd\circ \langle \operatorname{inv}^{\rhd}, \operatorname{inv} \rangle \big)_R(x),
    \end{aligned}
    \]
which implies 
\begin{equation}
\label{eq inv rhd}
    \operatorname{inv}^{\rhd}=\rhd\circ \langle \operatorname{inv}^{\rhd}, \operatorname{inv} \rangle.
\end{equation}
\end{remark}

\subsection{Coordinate rings of affine post-groups: co-post-Hopf algebras}
\label{subsection coordinate ring of affine post-group}

As the coordinate ring of an affine group is a commutative Hopf algebra, it is natural to ask what the corresponding algebraic structure is for affine post-groups. We use the Sweedler notation from Section \ref{sec affine algebraic group scheme}. 

\begin{definition}
\label{Def copostHopf algebra}
A \textbf{co-post-Hopf algebra} $(H, \blacktriangleright)$ over $k$ is a commutative Hopf $k$-algebra \((H, \mu, \eta, \Delta, \varepsilon, S)\) equipped with an algebra homomorphism 
\[
\blacktriangleright\colon  H \to H \otimes H, \quad h \mapsto h_{1_{\blacktriangleright}} \otimes h_{2_{\blacktriangleright}}
\]
that satisfies the following axioms:

\begin{enumerate}[label=(CPH\arabic*), ref=(CPH\arabic*)]
    \item \label{Item:CPH1} The following identity
    \begin{equation}
    \label{eq cph1}
        (\operatorname{id} \otimes \Delta) \circ \blacktriangleright = (\mu \otimes \operatorname{id} \otimes \operatorname{id}) \circ (\operatorname{id} \otimes \tau \otimes \operatorname{id}) \circ (\blacktriangleright \otimes \blacktriangleright) \circ \Delta 
    \end{equation}
    holds, where \(\tau(x \otimes y) = y \otimes x\). Equivalently, for any $h\in H$,
     \begin{equation*}
        h_{1_\blacktriangleright} \otimes h_{2_\blacktriangleright 1} \otimes h_{2_\blacktriangleright 2}=h_{1 1_\blacktriangleright} h_{2 1_\blacktriangleright} \otimes h_{1 2_\blacktriangleright} \otimes h_{2 2_\blacktriangleright}.
     \end{equation*}

    \item \label{Item:CPH2} For any \(R \in \mathsf{Alg}_k\) and \(\varphi \in \operatorname{Hom}_{\mathsf{Alg}_k}(H,R)\), the map
    \[
    L^{\blacktriangleright}_{\varphi}\colon  \operatorname{Hom}_{\mathsf{Alg}_k}(H, R) \to \operatorname{Hom}_{\mathsf{Alg}_k}(H, R), \qquad \psi \mapsto \mu_R \circ (\varphi \otimes \psi) \circ \blacktriangleright
    \]
    is bijective, where \(\mu_R\) denotes the multiplication of \(R\).

    \item \label{Item:CPH3} The following identity
    \begin{equation}
    \label{eq cph3}
        (\operatorname{id} \otimes \blacktriangleright) \circ \blacktriangleright = (\mu \otimes \operatorname{id} \otimes \operatorname{id}) \circ (\operatorname{id} \otimes \blacktriangleright \otimes \operatorname{id}) \circ (\Delta \otimes \operatorname{id}) \circ \blacktriangleright
    \end{equation}
    holds. Equivalently, for any $h \in H$, 
     \begin{equation*}
       h_{1_\blacktriangleright} \otimes h_{2_\blacktriangleright 1_\blacktriangleright} \otimes h_{2_\blacktriangleright 2_\blacktriangleright}=h_{1_\blacktriangleright 1} h_{1_\blacktriangleright 2 1_\blacktriangleright} \otimes h_{1_\blacktriangleright 2 2_\blacktriangleright} \otimes h_{2_\blacktriangleright}.
      \end{equation*}
\end{enumerate}
\end{definition}

\begin{remark}
As an example of our notation, $h_1 h_{2 1_\blacktriangleright} \otimes h_{2 2_\blacktriangleright}$ means first applying $\Delta$, then applying $\blacktriangleright$ to its second factor, and finally multiplying the first two factors. 
The above definition amounts to regarding a post-group object as living in $\mathsf{Alg}_k^{\operatorname{opp}}$, the opposite category of $\mathsf{Alg}_k$.
\end{remark}

\begin{definition}
A {\bf homomorphism} of co-post-Hopf algebras
\[
  \varphi\colon  (H, \blacktriangleright_H) \to (K,\blacktriangleright_K)
\]
is a Hopf algebra homomorphism \(\varphi\colon  H \to K\) satisfying $(\varphi \otimes \varphi) \circ \blacktriangleright_H = \blacktriangleright_K \circ \varphi$.
\end{definition}

We denote by \(\mathsf{CoPHopf}_k\) the category of co-post-Hopf algebras over \(k\).

\begin{example}
Let $(G,\cdot,\rhd)$ be a finite post-group. The algebra $k^G$ of $k$-valued functions on $G$ is a commutative Hopf algebra with pointwise multiplication and
\[
(\Delta f)(x,y)=f(x\cdot y),\qquad \varepsilon(f)=f(e),\qquad (Sf)(x)=f(x^{-1}).
\]
Under $k^G\otimes k^G\simeq k^{G\times G}$, the map
\[
\blacktriangleright\colon k^G\to k^G\otimes k^G,\qquad \blacktriangleright(f)(x,y)=f \circ \rhd (x,y)=f(x\rhd y)
\]
makes $(k^G,\blacktriangleright)$ a co-post-Hopf algebra.
\end{example}

\begin{proof}
First, the map $\blacktriangleright$ is an algebra homomorphism. Evaluating at $(x,y,z)$, conditions \ref{Item:CPH1} and \ref{Item:CPH3} follow from \ref{item:postgroup1} and \ref{item:postgroup3}, respectively.
For \ref{Item:CPH2}, fix $R\in\mathsf{Alg}_k$ and $\varphi\in\operatorname{Hom}_{\mathsf{Alg}_k}(k^G,R)$. Let $\delta_x$ be the characteristic function of $\{x\}$ and put $e_x=\varphi(\delta_x)$. Then $e_x^2=e_x$, $e_xe_y=0$ for $x\ne y$, and $\sum_xe_x=1$. Since
\[
   \blacktriangleright(f)=\sum_x \delta_x \otimes (f\circ L_x^\rhd),
\]
the definition of $L^{\blacktriangleright}_\varphi$ gives $L^{\blacktriangleright}_\varphi(\psi)(f) =\sum_x e_x\psi(f\circ L_x^\rhd)$. Define
\[
  T_\varphi(\psi)(f) = \sum_xe_x\psi\bigl(f\circ(L_x^\rhd)^{-1}\bigr).
\]
The map $T_\varphi(\psi)$ is $k$-linear and sends $1$ to $\sum_x e_x=1$. Since the $e_x$'s are orthogonal idempotents,
\begin{align*}
     T_\varphi(\psi)(f)\,T_\varphi(\psi)(g) &= \sum_x e_x \psi \bigl((f \circ (L_x^\rhd)^{-1}) (g \circ (L_x^\rhd)^{-1})\bigr)\\
     &= \sum_x e_x\psi\bigl((fg)\circ(L_x^\rhd)^{-1}\bigr) = T_\varphi(\psi)(fg).
\end{align*}
Thus $T_\varphi(\psi)$ is an algebra homomorphism. Moreover,
\[
\begin{aligned}
  L^{\blacktriangleright}_\varphi(T_\varphi(\psi))(f) &= \sum_{x,y}e_xe_y \psi\bigl(f\circ L_x^\rhd\circ(L_y^\rhd)^{-1}\bigr) = \psi(f),\\
  T_\varphi(L^{\blacktriangleright}_\varphi(\psi))(f) &= \sum_{x,y}e_xe_y \psi\bigl(f\circ(L_x^\rhd)^{-1}\circ L_y^\rhd\bigr) =\psi(f).
\end{aligned}
\]
Thus $T_\varphi$ is the inverse of $L^{\blacktriangleright}_\varphi$.
\end{proof}

\begin{proposition}
For any co-post-Hopf algebra, we have
\begin{equation}
\label{CPH prop 1}
\varepsilon(h_{1}) h_{2 1_\blacktriangleright} \otimes h_{2 2_\blacktriangleright} = h_{1_\blacktriangleright} \otimes h_{2_\blacktriangleright},
\end{equation}
where $h \in H$. Replacing $h$ by $h_{1_\blacktriangleright}$ gives
\begin{equation}
\label{CPH prop 1'}
\varepsilon(h_{1_\blacktriangleright 1}) \, h_{1_\blacktriangleright 2 1_\blacktriangleright} \otimes h_{1_\blacktriangleright 2 2_\blacktriangleright}
=
h_{1_\blacktriangleright 1_\blacktriangleright} \otimes h_{1_\blacktriangleright 2_\blacktriangleright}.
\end{equation}
\end{proposition}
\begin{proof}
By the \(k\)-linearity of \(\blacktriangleright\),
\[
\varepsilon(h_{1}) \, \blacktriangleright(h_{2}) = \blacktriangleright(\varepsilon(h_{1}) h_{2}).
\]
Using the counit identity \(\varepsilon(h_{1}) h_{2} = h\) on the right-hand side, we get
\[
\varepsilon(h_{1}) \, \blacktriangleright(h_{2}) = \blacktriangleright(h).
\]
Expanding both sides gives Eq.\eqref{CPH prop 1}. 
\end{proof}

Further properties of co-post-Hopf algebras will be discussed after establishing their equivalence with affine post-group schemes; see the next subsection.  

\begin{theorem}
\label{Thm affine post-group equivalent to co-post Hopf}
If \((G, m, \rhd)\) is an affine post-group, then  $\mathcal{O}(G, m, \rhd)$ is a co-post-Hopf algebra. Conversely, if $(H, \Delta, \blacktriangleright)$ is a co-post-Hopf algebra, then $(h^H,h^{\Delta},h^{\blacktriangleright})$ is an affine post-group. 
\end{theorem}

\begin{proof}
    Let $(G,m,\rhd)$ be an affine post-group. Since $G\times G\cong h^{\O(G)\otimes\O(G)}$, the natural transformations $m$, $e$, $\operatorname{inv}$, and $\tr$ correspond under the Yoneda lemma to $\Delta$, $\eps$, $S$, and a $k$-algebra homomorphism $\blacktriangleright\colon \O(G)\to\O(G)\otimes\O(G)$ satisfying
    \[
    \tr_R\langle f,g\rangle=\mu_R\circ(f\otimes g)\circ\blacktriangleright.
    \]
    By the Yoneda lemma, the axioms \ref{item postgroup object 1} and \ref{item postgroup object 3} correspond to the axioms \ref{Item:CPH1} and \ref{Item:CPH3} on the commutative Hopf algebra $(\O(G),\Delta,\eps,S)$.

    It remains to show the same for \ref{item postgroup object 2} and \ref{Item:CPH2}.
    By Proposition \ref{Prop: affine postgroup=postgroup object in repfunc}, $(G,m,\rhd)$ is a post-group object in $\mathsf{RepFunc}[\mathsf{Alg}_k, \mathsf{Set}]$. For every representable functor $F=h^A$, the Yoneda bijection $\iota\colon\Hom(\O(G),A)\to\Hom(F,G)$ sends $f\mapsto T_f$, where $(T_f)_R(g)=g\circ f$. Given $\varphi\in\Hom(\O(G),A)$, let $\Phi\coloneq\iota(\varphi)$. Then, for all $R\in\mathsf{Alg}_k$ and $g\in F(R)$, we have
    \[
\begin{aligned}
\big(\rhd \circ \langle \Phi, T_f \rangle \big)_R(g) 
&= \rhd_R \circ \big\langle\Phi_R, (T_f)_R\big\rangle(g)\\
&=\rhd_R \langle g \circ \varphi, g \circ f \rangle\\
&=\mu_R \circ \big( (g \circ \varphi)\otimes (g \circ f) \big) \circ \blacktriangleright\\
&=\mu_R \circ (g \otimes g)\circ (\varphi \otimes f) \circ \blacktriangleright\\
&= g\circ \mu_{A} \circ (\varphi \otimes f) \circ \blacktriangleright,
\end{aligned}
\]
which implies that
\[
\iota^{-1}\circ L^\tr_\Phi\circ\iota=L^\blacktriangleright_\varphi.
\]
Since $\iota$ is bijective, it follows that \ref{item postgroup object 2} and \ref{Item:CPH2} are equivalent.

The converse follows from reading the same correspondences backwards.
\end{proof}

\begin{theorem}
\label{Thm AffPGrp simeq CoHopf opp}
The assignments 
\[
\mathcal{O}\colon(G, m, \rhd) \mapsto \mathcal{O}(G, m, \rhd) \quad \text{and} \quad H \rightsquigarrow h^H \colon(H, \Delta,\blacktriangleright) \mapsto (h^H,h^{\Delta}, h^{\blacktriangleright})
\] 
given by Theorem \ref{Thm affine post-group equivalent to co-post Hopf} are quasi-inverse functors, establishing an equivalence of categories
\[
\mathsf{AffPGrp}_k \simeq \mathsf{CoPHopf}_k^{\mathrm{opp}}.
\]
\end{theorem}

\begin{proof}
    Recall that morphisms of affine groups $\Psi\colon G\to G'$ correspond bijectively to Hopf algebra homomorphisms $\psi\colon\O(G')\to\O(G)$. By the Yoneda lemma, the identities
    \[
    \Psi \circ \rhd=\rhd' \circ (\Psi \times \Psi),\qquad\blacktriangleright \circ \psi= (\psi \otimes \psi) \circ \blacktriangleright'
    \]
    are equivalent. Hence, the claim follows from Theorem \ref{Thm affine post-group equivalent to co-post Hopf}.
\end{proof}

\subsection{Co-post-Hopf algebras and their sub-adjacent Hopf algebras} 

We construct a new Hopf algebra from the co-post-Hopf algebra and further study its properties.

\begin{proposition}
\label{Prop sub adjacent Hopf algebra}
    Let $(H,\mu,\eta,\Delta,\varepsilon, S, 
    \blacktriangleright)$ be a co-post-Hopf algebra. Define 
    \begin{equation}
    \label{subjacent Hopf algebra Delta}
        \Delta_\blacktriangleright(h)=h_1 h_{2 1_\blacktriangleright} \otimes h_{2 2_\blacktriangleright},
    \end{equation}
    where $h \in H$.
    Then $(H,\mu,\eta, \Delta_\blacktriangleright, \varepsilon)$ is a commutative bialgebra. Moreover, there exists a unique antipode $S_\blacktriangleright\colon H \to H$ satisfying
    \begin{equation}
    \label{eq S blacktri}
        S_\blacktriangleright= \mu \circ (S_\blacktriangleright \otimes S) \circ \blacktriangleright
    \end{equation}
    such that $H_\blacktriangleright\coloneq (H,\mu,\eta, \Delta_\blacktriangleright, \varepsilon, S_\blacktriangleright)$ is a commutative Hopf algebra, called the \textbf{sub-adjacent Hopf algebra} of the co-post-Hopf algebra $(H,\mu,\eta,\Delta,\varepsilon, S, 
    \blacktriangleright)$.
\end{proposition}
\begin{proof}
    By  Theorem \ref{Thm AffPGrp simeq CoHopf opp}, $H$ corresponds to an affine post-group $(G,m,e,\operatorname{inv},\rhd)$. Let $(G,m^\rhd,e,\operatorname{inv}^\rhd)$ be the sub-adjacent affine group of $(G,m,e,\operatorname{inv},\rhd)$ by Proposition \ref{Prop sub-adjacent affine group of affine postgroup}, where 
    \[
    m^\rhd=m \circ (\operatorname{id} \times \rhd) \circ \big(\langle \operatorname{id}, \operatorname{id} \rangle \times \operatorname{id}\big)=h^\Delta \circ (\operatorname{id} \times h^\blacktriangleright) \circ \big(\langle \operatorname{id}, \operatorname{id} \rangle \times \operatorname{id}\big).
    \] 
    Then we have a commutative Hopf algebra $(H,\mu,\eta, \Delta_\blacktriangleright,\varepsilon, S_\blacktriangleright)\coloneq \mathcal{O}(G,m^\rhd,e,\operatorname{inv}^\rhd)$ as the coordinate ring of the affine group, where 
    \[
    \Delta_\blacktriangleright=\big(\big(\mu \circ(\operatorname{id}\otimes\operatorname{id})\big)\otimes \operatorname{id}\big) \circ (\operatorname{id} \otimes \blacktriangleright) \circ \Delta=(\mu \otimes \operatorname{id}) \circ (\operatorname{id} \otimes \blacktriangleright) \circ \Delta.
    \]
    Namely, for $h \in H$, $\Delta_\blacktriangleright(h)=h_1 h_{2 1_\blacktriangleright} \otimes h_{2 2_\blacktriangleright}$. Furthermore, the Yoneda lemma shows that there is a unique algebra homomorphism $S_\blacktriangleright\colon H \to H$ corresponding to the natural transformation $\operatorname{inv}^{\rhd}:G \to G$. Dualizing Eq. \eqref{eq inv rhd} in Remark~\ref{Remark on sub-adjacent affine group} yields Eq. \eqref{eq S blacktri}. The proof is complete.
\end{proof}

\begin{remark}
\label{Remark sub=sub}
    The proof of Proposition \ref{Prop sub adjacent Hopf algebra} implies the following diagram:
\begin{equation*}
\begin{tikzcd}
\text{affine post-group}
  \arrow[r, "\mathcal O(-)"]
  \arrow[d, "\text{sub-adjacent}"']
&
\substack{\text{co-post-Hopf algebra}}
  \arrow[d, "\text{sub-adjacent}"]
\\
\text{affine group}
  \arrow[r, "\mathcal O(-)"]
&
\substack{\text{commutative} \\ \text{Hopf algebra}}.
\end{tikzcd}
\end{equation*}
\end{remark}

\begin{example}
\label{ex:copost-hopf-basic}
Let $H=\mathcal O(G)$. By Theorem~\ref{Thm affine post-group equivalent to co-post Hopf}, the structures in
Example~\ref{ex:affine-post-basic} give the following
co-post-Hopf algebras $(H,\blacktriangleright)$.
Their sub-adjacent Hopf algebras are obtained from
Proposition~\ref{Prop sub adjacent Hopf algebra}.
\begin{enumerate}
    \item[(i)]
    The operation $\rhd\coloneq \pi_2: G \times G \to G$ corresponds to 
    \[
    \blacktriangleright\coloneq \iota_2\colon H \to H \otimes H, \quad h \mapsto 1\otimes h.
    \]
    Here $\Delta_{\blacktriangleright}=\Delta$ and
    $S_{\blacktriangleright}=S$, so
    $H_{\blacktriangleright}$ coincides with $H$.

    \item[(ii)]
    For $\rhd_R(f,g)=f^{-1}\cdot_R g\cdot_R f$, define
    \[
    \blacktriangleright(h)=S(h_1)h_3\otimes h_2.
    \]
    For every $R\in\mathsf{Alg}_k$ and $f,g\in G(R)$,
    \begin{align*}
      \rhd_R(f,g) &=\mu_R \circ (\mu_R \otimes \operatorname{id}) \circ \bigl((f \circ S) \otimes g \otimes f\bigr) \circ (\Delta \otimes \operatorname{id}) \circ \Delta\\
      &=\mu_R\circ(f\otimes g)\circ\blacktriangleright,
    \end{align*}
which implies, by the Yoneda lemma, $\blacktriangleright$ corresponds to $\rhd$. Here $\Delta_{\blacktriangleright}=\tau\circ\Delta$ and $S_{\blacktriangleright}=S$.
\end{enumerate}
\end{example}

We now establish other basic properties of co-post-Hopf algebras.

\begin{proposition}
For any co-post-Hopf algebra \((H, \mu, \eta, \Delta, \varepsilon, S, \blacktriangleright)\), the following identities hold for all \(h \in H\):
\begin{equation}
\label{CPH counit 1}
\varepsilon(h_{1_\blacktriangleright})h_{2_\blacktriangleright}=h ,
\end{equation}
\begin{equation}
\label{CPH counit 2}
    h_{1_\blacktriangleright} \varepsilon(h_{2_\blacktriangleright})=\varepsilon(h)1_H,
\end{equation}
\begin{equation}
\label{CPH counit 3}
\varepsilon(h_{1_\blacktriangleright}) \varepsilon(h_{2_\blacktriangleright})=\varepsilon(h),
\end{equation}
where \(\blacktriangleright(h)=h_{1_\blacktriangleright} \otimes h_{2_\blacktriangleright}\) and \(1_H\) is the identity element of \(H\).
\end{proposition}
\begin{proof}
Let $(G,m,e,\operatorname{inv},\rhd)$ be an affine post-group corresponding to the co-post-Hopf algebra $(H, \mu, \eta, \Delta, \varepsilon, S, \blacktriangleright)$ via Theorem \ref{Thm AffPGrp simeq CoHopf opp}. By dualizing the two diagrams in \eqref{Diag f rhd 1=1 and 1 rhd f = f}, we obtain the following 
\[
        \begin{tikzcd}[column sep=2cm]
            H \ar[equal,r] \ar[d, "\blacktriangleright"] & H \\
            H \otimes H \ar[r, "\varepsilon \otimes \operatorname{id}"] & k \otimes H \ar[u, "{\mu \circ (\eta \otimes \operatorname{id})}"],
        \end{tikzcd}
        \quad \text{and} \quad 
        \begin{tikzcd}[column sep=2cm]
            H \ar[d, "\blacktriangleright"] \ar[r, "\eta \circ \varepsilon"] & H \ar[equal, d] \\
            H \otimes H \ar[r, "{\mu \circ \big(\operatorname{id} \otimes (\eta \circ \varepsilon) \big)}"] & H
        \end{tikzcd}
\]
which are equivalent to Eqs. \eqref{CPH counit 1} and \eqref{CPH counit 2}, respectively. Moreover, Eq. \eqref{CPH counit 3} is obtained by applying $\varepsilon$ on Eq. \eqref{CPH counit 1}.
\end{proof}

\section{An equivalent description via affine skew braces} \label{sec affine post-group schemes and related affine schemes}

In \cite{AGV17}, Angiono, Galindo and Vendramin briefly observed that affine skew braces are related to the Yang--Baxter equation in the setting of affine schemes.  
We prove that the category of affine post-groups is isomorphic to that of affine skew braces, thereby linking affine post-groups to affine Yang--Baxter schemes.

\subsection{Post-groups and skew braces}
Recall from \cite[Definition 1.1]{GV17} that a set $A$ with two binary operations $\cdot$ and $\circ$ is called a \textbf{skew left brace} if $(A, \cdot)$ and $(A, \circ)$ are groups such that
\[
a \circ (b \cdot c)= (a \circ b) \cdot a^{-1} \cdot (a \circ c), \quad \forall a,b,c \in A.
\]
Here $a^{-1}$ denotes the inverse of $a$ with respect to $\cdot$ and we refer to a skew left brace simply as a \textbf{skew brace} in this paper. A skew brace $(A,\cdot, \circ)$ is called a \textbf{brace} if the group $(A, \cdot)$ is abelian.

\begin{proposition}[\cite{BGST24}]
\label{Prop postgroup-skewbrace}
The category of post-groups is isomorphic to the category of skew braces. More precisely:
\begin{enumerate}
    \item[\rm (1)] Given a post-group \((G, \cdot, \rhd)\), the operation \(a \circ b = a \cdot (a \rhd b)\) defines a skew brace \((G, \cdot, \circ)\).
    \item[\rm (2)] Conversely, given a skew brace \((G, \cdot, \circ)\), the operation \(a \rhd b = a^{-1} \cdot (a \circ b)\) defines a post-group \((G, \cdot, \rhd)\).
    \item[\rm (3)] These constructions are inverse to each other and functorial, yielding an isomorphism of categories \(\mathsf{PGrp} \simeq \mathsf{SBr}\), where $\mathsf{SBr}$ denotes the category of skew braces.
\end{enumerate}
\end{proposition}
\begin{proof}
    See \cite[Proposition 3.22, Proposition 3.24 and Theorem 3.25]{BGST24}.
\end{proof}

\subsection{Affine post-groups, affine skew braces and their coordinate rings}
We begin by recalling the notion of affine skew braces. 

\begin{definition}{\rm (\cite{AGV17})}
\label{def affine skew brace}
    An \textbf{affine (skew) brace scheme over $k$} (or simply an \textbf{affine (skew) brace}) is a representable functor $\mathcal{A} \colon \mathsf{Alg}_k \to \mathsf{Set}$ together with two natural transformations \(m_\bullet, m_\circ : \mathcal{A} \times \mathcal{A} \to \mathcal{A}\) such that, for all $R$ in $\mathsf{Alg}_k$,
\[
{m_\bullet}_R : \mathcal{A}(R) \times \mathcal{A}(R) \to \mathcal{A}(R),\quad (a, b) \mapsto a \bullet_R b
\]
and
\[
{m_\circ}_R : \mathcal{A}(R) \times \mathcal{A}(R) \to \mathcal{A}(R),\quad (a, b) \mapsto a \circ_R b
\]
   make \(\big(\mathcal{A}(R), \bullet_R, \circ_R \big)\) into a (skew) brace. If $\mathcal{A}$ is represented by a finitely generated $k$-algebra, then it is called an \textbf{affine algebraic (skew) brace}.
\end{definition}

\begin{remark}
Definition \ref{def affine skew brace} may be regarded as a finer version of this terminology, in that it distinguishes affine skew braces from affine braces. More precisely, our affine skew brace agrees with the original notion of an affine brace introduced by Angiono, Galindo and Vendramin \cite[Definition 5.1]{AGV17}, while our affine brace is obtained by further requiring commutativity. Moreover, an affine (skew) brace is just a (skew) brace object in $\mathsf{RepFunc}[\mathsf{Alg}_k, \mathsf{Set}]$; see \cite{MZL26}.
\end{remark}

A \textbf{homomorphism of affine (skew) brace} \(\Phi\colon  (\mathcal{A}, m_\bullet, m_\circ) \to (\mathcal{A}', m_\bullet', m_\circ')\) is a natural transformation such that, for any $R \in \mathsf{Alg}_k$, the map
\[
\Phi_R\colon  \mathcal{A}(R) \to \mathcal{A}'(R)
\]
is a homomorphism of (skew) braces. Thus, affine (skew) braces over \(k\) together with their homomorphisms form a category, denoted by \(\mathsf{AffSBr}_k\) (resp. \(\mathsf{AffBr}_k\)).

\begin{proposition}
\label{Prop affine postgroup give affine skew brace}
Let \((G, m, \rhd)\) be an affine post-group and \((G, m^\rhd)\) be its sub-adjacent affine group. Then \((G, m, m^\rhd)\) is an affine skew brace. Moreover, let $\Psi\colon (G,m,\rhd) \to (G',m',\rhd')$ be a homomorphism of affine post-groups. Then $\Psi$ induces a homomorphism of affine skew braces from $(G,m,m^{\rhd})$ to $(G',m',m^{\rhd'})$.
\end{proposition}
\begin{proof}
First, for each \(R \in \mathsf{Alg}_k\), the triple \(\big(G(R), m_R, \rhd_R \big)\) is a post-group and the pair$\big(G(R),m^{\rhd}_R \big)$ is a group, where $m^{\rhd}_R(f,g)=f \circ_R g=f \cdot_R (f \rhd_R g)$. One can check that
\[
f \circ_R(g\cdot_Rh)=(f\circ_Rg) \cdot_R f^{-1} \cdot_R(f\circ_R h),
\]
which implies $\big ( G(R), \cdot_R, \circ_R \big)$ is a skew brace.
Therefore, \((G, m, m^\rhd)\) is an affine skew brace. Next, for each homomorphism of post-groups $\Psi_R$ from $\big(G(R),m_R,\rhd_R\big)$ to $\big(G'(R),m'_R,\rhd'_R\big)$, we have
\begin{align*}
    \Psi_R \circ m^{\rhd}_R(x,y)&=\Psi_R (x \cdot_R (x \rhd_R y))\\
    &=\Psi_R(x) \cdot'_R \big( \Psi_R(x) \rhd'_R \Psi_R(y) \big)\\
    &=m^{\rhd'}_R \big( \Psi_R(x) , \Psi_R(y) \big),
\end{align*}
which implies $\Psi_R$ is also a homomorphism of skew braces from $\big(G(R),m_R,m^{\rhd}_R\big)$ to $\big(G'(R),m'_R,m^{\rhd'}_R\big)$. That is, $\Psi$ is a homomorphism of affine skew brace from $(G,m,m^{\rhd})$ to $(G',m',m^{\rhd'})$.
\end{proof}

\begin{proposition}
\label{Prop affine skew brace give affine postgroup}
Let $(\mathcal{A},m_\bullet,m_\circ)$ be an affine skew brace. Define a natural transformation $\rhd:\mathcal{A} \times \mathcal{A} \to \mathcal{A}$ by the following commutative diagram:
\begin{equation}
\label{rhd from affine skew brace diagram}
    \begin{tikzcd}[column sep=2cm]
       \mathcal{A}^{\times 2} \ar[r, "\rhd"] \ar[d,"{\langle \operatorname{id},\operatorname{id}\rangle \times \operatorname{id}}"] & \mathcal{A} \\
       \mathcal{A}^{\times 3} \ar[r, "{\operatorname{inv}_\bullet \times m_\circ}"] & \mathcal{A}^{\times 2} \ar[u, "{m_\bullet}"].
    \end{tikzcd}
\end{equation}
Then $(\mathcal{A}, m_\bullet, \rhd)$ is an affine post-group. Moreover, let $\Psi\colon (\mathcal{A},m_\bullet,m_\circ) \to (\mathcal{A}',m'_\bullet,m'_\circ)$ be a homomorphism of affine skew braces. Then $\Psi$ is also a homomorphism of affine post-groups from $(\mathcal{A}, m_\bullet, \rhd)$ to $(\mathcal{A}', m'_\bullet, \rhd')$.
\end{proposition}
\begin{proof}
The diagram \eqref{rhd from affine skew brace diagram} means that, for each $R \in \mathsf{Alg}_k$ and $x,y \in \mathcal{A}(R)$,
\[
\rhd_R(x, y) = x^{-1} \cdot_R (x \circ_R y),
\]
which implies $\big(\mathcal{A}(R), \cdot_R, \rhd_R\big)$ is a post-group by Proposition \ref{Prop postgroup-skewbrace}. Thus, $(\mathcal{A}, m_\bullet, \rhd)$ is an affine post-group. And, for each $\Psi_R$, we have
\begin{align*}
    \Psi_R \circ \rhd_R(x,y)&=\Psi_R\big(x^{-1} \cdot_R (x \circ_R y)\big)\\
    &=\Psi_R(x)^{-1} \cdot'_R \big( \Psi_R(x) \circ'_R \Psi_R(y) \big)\\
    &=\Psi_R(x) \rhd_R' \Psi_R(y),
\end{align*}
which implies $\Psi_R$ is also a homomorphism of post-groups from $(\mathcal{A}(R), m_{\bullet R}, \rhd_R)$ to $(\mathcal{A}'(R), m'_{\bullet R}, \rhd'_R)$. Namely, $\Psi$ is a homomorphism of affine post-groups from $(\mathcal{A}, m_\bullet, \rhd)$ to $(\mathcal{A}', m'_\bullet, \rhd')$.
\end{proof}

\begin{theorem}
\label{Thm AffPGrp=AffSBr}
Proposition \ref{Prop affine postgroup give affine skew brace} gives a functor $\mathcal{F}: \mathsf{AffPGrp}_k \to \mathsf{AffSBr}_k$, and Proposition \ref{Prop affine skew brace give affine postgroup} gives a functor $\mathcal{P}\colon \mathsf{AffSBr}_k \to \mathsf{AffPGrp}_k$. Moreover, they give an isomorphism between the categories $\mathsf{AffPGrp}_k$ and $\mathsf{AffSBr}_k$.
\end{theorem}

\begin{proof}
First, the functors $\mathcal{F}$ and $\mathcal{P}$ are well-defined. Let $(G,m,\rhd)$ and $(G',m',\rhd')$ be two affine post-groups, and $\Psi$ be a homomorphism of affine post-groups from $(G,m,\rhd)$ to $(G',m',\rhd')$. 
By Propositions \ref{Prop affine postgroup give affine skew brace} and \ref{Prop affine skew brace give affine postgroup}, we have $\mathcal{P} \mathcal{F}(G,m,\rhd)=(G,m,\rhd)$ and $\mathcal{P} \mathcal{F}(\Psi)=\Psi$. So $\mathcal{P} \mathcal{F}= \operatorname{id}$. Similarly, $\mathcal{F}\mathcal{P}= \operatorname{id}$.
Thus, we obtain an isomorphism $\mathsf{AffPGrp}_k \simeq \mathsf{AffSBr}_k$.
\end{proof}

We next give a brief description of Hopf co-braces without recalling the precise definition; for details we refer to \cite{AGV17}. We denote by \(\mathsf{HopfCoBr}_k\) the category of Hopf co-braces. In fact, Hopf co-braces are just coordinate rings of affine skew braces. Thus, we have the following equivalence between categories \cite[Proposition 5.6]{AGV17}:
\[
\mathsf{AffSBr}_k \simeq \mathsf{HopfCoBr}_k^{\mathrm{opp}}.
\]
Combining this with Theorems~\ref{Thm AffPGrp simeq CoHopf opp} and Theorem \ref{Thm AffPGrp=AffSBr}, we obtain the following chain of equivalences
\[
\mathsf{CoPHopf}_k \;\simeq\; \mathsf{AffPGrp}_k^{\mathrm{opp}} \;\simeq\; \mathsf{AffSBr}_k^{\mathrm{opp}} \;\simeq\; \mathsf{HopfCoBr}_k,
\]
which implies an equivalence between co-post-Hopf algebras and Hopf co-braces. It follows directly that a co-post-Hopf algebra, together with its sub-adjacent Hopf algebra, gives rise to a Hopf co-brace.

\begin{remark}
Recently, Li, Sheng and Tang \cite{LST24} introduced (cocommutative) \textbf{post-Hopf algebras}, which are equivalent to cocommutative Hopf braces. Moreover, Hopf co-braces are dual to Hopf braces \cite{AGV17}. Hence, a co-post-Hopf algebra can be naturally understood as a dualization of a post-Hopf algebra.
\end{remark}

\subsection{Application to affine Yang--Baxter schemes}

Every affine post-group $(G,m,\rhd)$ induces a nondegenerate \textbf{affine Yang--Baxter scheme} (see \cite{MZL26,Ta25} for the definition of nondegenerate affine Yang--Baxter schemes). 
In more detail, given $R$ in $\mathsf{Alg}_k$ and $a,b \in G(R)$, define
    \begin{equation*}
        \mathcal{K}_R(a,b)=\big(a \rhd_R b, (L^{\rhd_R}_{a \rhd_R b})^{-1} \big((a \rhd_R b)^{-1} \cdot_R a  \cdot_R (a \rhd_R b) \big)\big).
    \end{equation*}
    Then $\mathcal{K}_R$ is natural in $R$. Moreover, $(G,\mathcal{K})$ is a nondegenerate affine Yang--Baxter scheme, and it is involutive if and only if $(G,m)$ is an abelian affine group.
    This follows from Theorem
    \ref{Thm AffPGrp=AffSBr} and \cite[Prop.\ 5.22]{MZL26} after translating the operations of $(G,m,\rhd)$ into those in an affine skew brace via
 \[
 a \cdot_R (a \rhd_R b)= a \circ_R b = a \cdot _R \lambda_a(b),\qquad\lambda_a(b)\coloneq a^{-1}\cdot(a\circ b).
 \]

\section{Post-Lie algebras and Rota--Baxter Lie algebras arising from affine group schemes}
\label{sec post-Lie and RB Lie}
In this section, we construct Lie-algebraic structures associated to two classes of affine group schemes with additional operations: affine post-group schemes and affine Rota--Baxter group schemes. More precisely, we show that the tangent space at the identity of an affine post-group scheme carries a natural post-Lie algebra structure, while that of an affine Rota--Baxter group scheme admits a natural Rota--Baxter Lie algebra structure. 
As a preparation, we first recall the definitions of post-Lie algebras and Rota--Baxter Lie algebras, as well as the Lie algebra of an affine group scheme.

\subsection{Post-Lie algebras and Rota--Baxter Lie algebras}

Recall \cite{V07} that a {\bf post-Lie algebra} consists of a Lie algebra $(\mathfrak g, [\cdot, \cdot]_{\mathfrak{g}})$ together with a binary operation $\triangleright\colon  \mathfrak g \otimes \mathfrak g \to \mathfrak g$ satisfying
\begin{align*}
x \triangleright [y, z]_{\mathfrak{g}} &= [x \triangleright y, z]_{\mathfrak{g}} + [y, x \triangleright z]_{\mathfrak{g}},  \\
([x, y]_{\mathfrak{g}} + x \triangleright y - y \triangleright x) \triangleright z
&= x \triangleright (y \triangleright z) - y \triangleright (x \triangleright z)
\end{align*}
for all $x, y, z \in \mathfrak g$. 
When $(\mathfrak{g}, [\cdot, \cdot]_{\mathfrak{g}})$ is abelian, this reduces to the notion of a \textbf{pre-Lie algebra}. 
Moreover, any post-Lie algebra gives rise to a new Lie algebra $\mathfrak g_{\triangleright}$ with bracket
\[
[x, y]_{\mathfrak{g}_{\triangleright}} \coloneq  [x, y]_{\mathfrak{g}} + x \triangleright y - y \triangleright x,
\] 
called, following \cite{BGST24,L26}, the \textbf{sub-adjacent Lie algebra} of $(\g, [\cdot,\cdot]_{\g}, \triangleright)$ (also known as the \textbf{Grossman--Larson Lie algebra} \cite{AEM24}). The bracket $[\cdot, \cdot]_{\g_{\triangleright}}$ is also called the post-Lie bracket \cite{BK25} or the composition Lie bracket
\cite{JZ25}.

A \textbf{homomorphism} $\phi$ from a post-Lie algebra $(\mathfrak{g}, [\cdot, \cdot], \triangleright)$ to a post-Lie algebra $(\mathfrak{g}',[\cdot, \cdot]', \triangleright')$ is a Lie algebra morphism $\phi\colon  \mathfrak{g} \to \mathfrak{g}'$ such that $\phi (x \triangleright y)= \phi(x) \triangleright' \phi(y)$ for all $x,y \in \mathfrak{g}$.

A closely related algebraic structure is that of Rota--Baxter Lie algebras; their relation to post-Lie algebras was established in\cite{BGN10}.
A linear map $\mathscr{B}: \mathfrak{g} \to \mathfrak{g}$ is called a \textbf{Rota--Baxter operator of weight $\lambda$} on a Lie algebra $\mathfrak{g}$ if 
\begin{equation*}
    [\mathscr{B}(x),\mathscr{B}(y)]_{\mathfrak{g}}= \mathscr{B}([\mathscr{B}(x),y]_{\mathfrak{g}} + [x,\mathscr{B}(y)]_{\mathfrak{g}} + \lambda[x,y]_{\mathfrak{g}}).
\end{equation*}
We call $(\mathfrak{g}, [\cdot, \cdot]_{\mathfrak{g}}, \mathscr{B})$ the \textbf{Rota--Baxter Lie algebra} of weight $\lambda$. Moreover, there is a new Lie bracket $[\cdot, \cdot]_{\mathscr{B}}$ on $\mathfrak{g}$ defined by
\[
[x, y]_{\mathscr{B}} = [\mathscr{B}(x),y]_{\mathfrak{g}} + [x,\mathscr{B}(y)]_{\mathfrak{g}} + \lambda[x,y]_{\mathfrak{g}}.
\]
Following \cite{GLS21}, the Lie algebra $(\mathfrak{g},[\cdot, \cdot]_{\mathscr{B}})$ is called the \textbf{descendant Lie algebra} of $(\mathfrak{g}, [\cdot, \cdot]_{\mathfrak{g}}, \mathscr{B})$, denoted by $\mathfrak{g}_{\mathscr{B}}$.

A homomorphism $\psi$ from a Rota--Baxter Lie algebra $(\mathfrak{g}, [\cdot, \cdot]_\mathfrak{g}, \mathscr{B})$ to a Rota--Baxter Lie algebra $(\mathfrak{h}, [\cdot, \cdot]_\mathfrak{h}, \mathscr{R})$ is a Lie algebra homomorphism $\psi\colon  \mathfrak{g} \to \mathfrak{h}$ such that $\psi \circ \mathscr{B}= \mathscr{R} \circ \psi$.

\begin{remark}
\label{Remark cB is also RB operator}
Let \(\mathscr B\) be a Rota--Baxter operator of weight \(\lambda\) on a Lie algebra \(\mathfrak g\). Then for any \(c \in k\), the operator \(c\mathscr B\) is again a Rota--Baxter operator on \(\mathfrak g\), but of weight \(c\lambda\). 
\end{remark}

\begin{remark}\label{rmk:splitting}
A Rota--Baxter Lie algebra (of weight $\lambda$) \((\mathfrak{g}, [\cdot, \cdot]_{\mathfrak{g}}, \mathscr{B})\) naturally induces a post-Lie algebra \((\mathfrak{g}, [\cdot, \cdot]_\lambda, \triangleright)\) by
\[
[\cdot, \cdot]_{\lambda}\coloneq  \lambda[\cdot, \cdot]_{\mathfrak{g}}, \quad x \triangleright y \coloneq  [\mathscr{B}(x), y]_{\mathfrak{g}}.
\]
In the case of $\lambda=1$, the induced post-Lie algebra $(\mathfrak{g}, [\cdot, \cdot]_1, \triangleright)$, where $[\cdot, \cdot]_1=[\cdot, \cdot]_{\mathfrak{g}}$, is called the splitting post-Lie algebra of \((\mathfrak{g}, [\cdot, \cdot]_{\mathfrak{g}}, \mathscr{B})\) (see \cite[Corollary 5.6]{BGN10} for the construction and \cite[Section 4.2]{GLS21} for the splitting post-Lie algebra).
\end{remark}

\subsection{The post-Lie algebra of an affine post-group scheme}

Recall (from, for example, \cite{W1979} or \cite{Milne17}) that the Lie algebra $\mathfrak{g}=\Lie(G)$ of an affine group scheme $G$ can be defined to be the vector space of $k$-derivations $\Der_{k,\eps}(\O(G),k)$---that is, the vector space of $k$-linear maps $d\colon\O(G)\to k$ satisfying
\[
d(xy)=d(x)\eps(y)+\eps(x)d(y)
\]
for all $x,y\in\O(G)$---equipped with the Lie bracket
\begin{equation}
\label{Lie bracket of k-derivation}
    [d_1,d_2](x)\coloneq d_1(x_1)d_2(x_2)-d_2(x_1)d_1(x_2),\qquad d_1,d_2\in\Lie(G).
\end{equation}

Now, let $(G, m, \rhd)$ be an affine post-group over $k$, and let
\[
\mathcal{O}(G, m, \rhd) = (H, \mu, \eta, \Delta, \varepsilon, S, \blacktriangleright)=(H, \blacktriangleright)
\]
be its coordinate ring. 
For $d_1, d_2 \in \mathfrak g$, define a binary operation 
\[
\triangleright\colon  \mathfrak{g} \otimes \mathfrak{g} \to \mathfrak{g}, \quad d_1 \otimes d_2 \mapsto \mu_k \circ (d_1 \otimes d_2)\circ \blacktriangleright.
\]
Namely, for $h \in H$,
\begin{equation}
\label{triangleright of k-derivation}
(d_1 \triangleright d_2)(h)
= d_1(h_{1_\blacktriangleright}) \, d_2(h_{2_\blacktriangleright}),
\end{equation}
where $\blacktriangleright(h) = h_{1_\blacktriangleright} \otimes h_{2_\blacktriangleright}$.
The bilinearity of $\triangleright\colon  \mathfrak{g} \times \mathfrak{g} \to \mathfrak{g}$ is evident. We now verify that $d_1 \triangleright d_2$ is again a $k$-derivation.

\begin{proposition}\label{prop:d1d2g}
For any $d_1, d_2 \in \mathfrak g$, we have
$d_1 \triangleright d_2 \in \mathfrak g$.
\end{proposition}

\begin{proof}
First, $d_1 \triangleright d_2$ is a $k$-linear map. Then, for any $x, y \in H$, using that $\blacktriangleright$ is an algebra homomorphism,
\[
(d_1 \triangleright d_2)(xy) = d_1\big((xy)_{1_\blacktriangleright}\big) \, d_2\big((xy)_{2_\blacktriangleright}\big) = d_1(x_{1_\blacktriangleright} y_{1_\blacktriangleright}) \,
   d_2(x_{2_\blacktriangleright} y_{2_\blacktriangleright}).
\]
Since $d_1$ and $d_2$ are $k$-derivations, we have
\begin{align*}
& \quad (d_1 \triangleright d_2)(xy)\\
&=\big(d_1(x_{1_\blacktriangleright})\varepsilon(y_{1_\blacktriangleright}) + \varepsilon(x_{1_\blacktriangleright}) d_1(y_{1_\blacktriangleright}) \big) \big(d_2(x_{2_\blacktriangleright})\varepsilon(y_{2_\blacktriangleright}) + \varepsilon(x_{2_\blacktriangleright}) d_2(y_{2_\blacktriangleright}) \big)\\
&=\underline{d_1(x_{1_\blacktriangleright}) d_2(x_{2_\blacktriangleright})} \,
\uuline{\varepsilon(y_{1_\blacktriangleright}) \varepsilon(y_{2_\blacktriangleright})} + \uwave{d_1(x_{1_\blacktriangleright}) \varepsilon(x_{2_\blacktriangleright})} \,
\underbrace{\varepsilon(y_{1_\blacktriangleright}) d_2(y_{2_\blacktriangleright})} \\
&\quad + \underbrace{\varepsilon(x_{1_\blacktriangleright}) d_2(x_{2_\blacktriangleright})} \, \uwave{d_1(y_{1_\blacktriangleright}) \,\varepsilon(y_{2_\blacktriangleright})} + \uuline{\varepsilon(x_{1_\blacktriangleright}) \varepsilon(x_{2_\blacktriangleright})} \,
\underline{d_1(y_{1_\blacktriangleright}) d_2(y_{2_\blacktriangleright})} \\
&= (d_1 \triangleright d_2)(x) \, \varepsilon(y) +\varepsilon(x) d_1(1_H) d_2(y) + d_2(x) \varepsilon(y) d_1(1_H) + \varepsilon(x) \, (d_1 \triangleright d_2)(y) \tag*{by Eqs.\eqref{CPH counit 1}, \eqref{CPH counit 2}, \eqref{CPH counit 3} and \eqref{triangleright of k-derivation}} \\
&= (d_1 \triangleright d_2)(x) \, \varepsilon(y)
   + \varepsilon(x) \, (d_1 \triangleright d_2)(y).\tag*{by $d_1(1_H)=0$}
\end{align*}
Thus \(d_1 \triangleright d_2 \in \operatorname{Der}_{k,\varepsilon}(H, k)=\mathfrak{g}\).
\end{proof}

\begin{proposition}
The operation $\triangleright$ satisfies the derivation compatibility
with the Lie bracket:
\[
d_1 \triangleright [d_2, d_3]
=
[d_1 \triangleright d_2, d_3]
+
[d_2, d_1 \triangleright d_3].
\]
\end{proposition}

\begin{proof}
Let \(h \in H\). We first compute the left-hand side:
\begin{align*}
 & \quad (d_{1}\triangleright [d_{2},d_{3}])(h)\\
 &= d_{1}(h_{1_{\blacktriangleright}})\,[d_{2},d_{3}](h_{2_{\blacktriangleright}})\\
 &=\underline{d_{1}(h_{1_{\blacktriangleright}}) d_{2}(h_{2_{\blacktriangleright}1})\,d_{3}(h_{2_{\blacktriangleright}2})}
- \underline{d_{1}(h_{1_{\blacktriangleright}}) d_{3}(h_{2_{\blacktriangleright}1})\,d_{2}(h_{2_{\blacktriangleright}2})}\\
&=d_1(h_{11_{\blacktriangleright}}\,h_{21_{\blacktriangleright}}) d_2(h_{12_{\blacktriangleright}}) d_3(h_{22_{\blacktriangleright}}) - d_1(h_{11_{\blacktriangleright}}\,h_{21_{\blacktriangleright}}) d_3(h_{12_{\blacktriangleright}}) d_2(h_{22_{\blacktriangleright}}) \tag*{by \ref{Item:CPH1}} \\
 & = d_1(h_{11_{\blacktriangleright}}) \underline{\varepsilon(h_{21_{\blacktriangleright}})} d_2(h_{12_{\blacktriangleright}}) \underline{d_3(h_{22_{\blacktriangleright}})} + \underline{\varepsilon(h_{11_{\blacktriangleright}})} d_1(h_{21_{\blacktriangleright}}) \underline{d_2(h_{12_{\blacktriangleright}})} d_3(h_{22_{\blacktriangleright}})\\
 & \quad - d_1(h_{11_{\blacktriangleright}}) \underline{\varepsilon(h_{21_{\blacktriangleright}})} d_3(h_{12_{\blacktriangleright}}) \underline{d_2(h_{22_{\blacktriangleright}})} - \underline{\varepsilon(h_{11_{\blacktriangleright}})} d_1(h_{21_{\blacktriangleright}}) \underline{d_3(h_{12_{\blacktriangleright}})} d_2(h_{22_{\blacktriangleright}})\\
 & = d_1(h_{11_{\blacktriangleright}}) d_2(h_{12_{\blacktriangleright}}) d_3(h_2) + d_1(h_{21_{\blacktriangleright}}) d_2(h_1) d_3(h_{22_{\blacktriangleright}})\\
 & \quad - d_1(h_{11_{\blacktriangleright}}) d_3(h_{12_{\blacktriangleright}}) d_2(h_2) - d_1(h_{21_{\blacktriangleright}}) d_3(h_1) d_2(h_{22_{\blacktriangleright}}). \tag*{by Eq.\eqref{CPH counit 1}}
\end{align*}
On the other hand, we have
\begin{align*}
[d_{1}\triangleright d_{2},d_{3}](h)
&= (d_{1}\triangleright d_{2})(h_{1})\,d_{3}(h_{2})
- d_{3}(h_{1})\,(d_{1}\triangleright d_{2})(h_{2})\\
&=d_{1}(h_{11_{\blacktriangleright}})\,d_{2}(h_{12_{\blacktriangleright}})\,d_{3}(h_{2}) - d_{3}(h_{1})\,d_{1}(h_{21_{\blacktriangleright}})\,d_{2}(h_{22_{\blacktriangleright}}).
\end{align*}
Similarly, 
\[
[d_{2},d_{1}\triangleright d_{3}](h)
= d_{2}(h_{1})\,d_{1}(h_{21_{\blacktriangleright}})\,d_{3}(h_{22_{\blacktriangleright}})
- d_{1}(h_{11_{\blacktriangleright}})\,d_{3}(h_{12_{\blacktriangleright}})\,d_{2}(h_{2}).
\]
Hence, we obtain
\[
d_{1}\triangleright [d_{2},d_{3}]
=
[d_{1}\triangleright d_{2},d_{3}]
+
[d_{2},d_{1}\triangleright d_{3}],
\]
as desired.
\end{proof}

\begin{proposition}
The operation $\triangleright$ satisfies the following identity:
\[
([d_1, d_2] + d_1 \triangleright d_2 - d_2 \triangleright d_1) \triangleright d_3
=
d_1 \triangleright (d_2 \triangleright d_3)
-
d_2 \triangleright (d_1 \triangleright d_3),
\]
where $d_1,d_2,d_3 \in \mathfrak{g}$.
\end{proposition}

\begin{proof}
Let $h \in H$. First, by definitions, expanding the left-hand side gives
\[
\begin{aligned}
\text{LHS}(h)
&= \big( [d_1, d_2](h_{1_\blacktriangleright}) + (d_1 \triangleright d_2) (h_{1_\blacktriangleright}) - (d_2 \triangleright d_1) (h_{1_\blacktriangleright}) \big) d_3(h_{2_\blacktriangleright}) \\
&=d_{1}(h_{1_{\blacktriangleright}1})\,d_{2}(h_{1_{\blacktriangleright}2})\,d_{3}(h_{2_{\blacktriangleright}})
-d_{2}(h_{1_{\blacktriangleright}1})\,d_{1}(h_{1_{\blacktriangleright}2})\,d_{3}(h_{2_{\blacktriangleright}})\\
&\quad +d_{1}(h_{1_{\blacktriangleright}1_{\blacktriangleright}})\,d_{2}(h_{1_{\blacktriangleright}2_{\blacktriangleright}})\,d_{3}(h_{2_{\blacktriangleright}})
-d_{2}(h_{1_{\blacktriangleright}1_{\blacktriangleright}})\,d_{1}(h_{1_{\blacktriangleright}2_{\blacktriangleright}})\,d_{3}(h_{2_{\blacktriangleright}}).
\end{aligned}
\]
And the right-hand side expands to
\begin{align*}
\text{RHS}(h)
&= d_{1}(h_{1_{\blacktriangleright}})\,d_{2}(h_{2_{\blacktriangleright}1_{\blacktriangleright}})\,d_{3}(h_{2_{\blacktriangleright}2_{\blacktriangleright}}) - d_{2}(h_{1_{\blacktriangleright}})\,d_{1}(h_{2_{\blacktriangleright}1_{\blacktriangleright}})\,d_{3}(h_{2_{\blacktriangleright}2_{\blacktriangleright}}).
\end{align*}
Then, we compute the first summand of $\text{RHS}(h)$:
\begin{align*}
& \quad d_{1}(h_{1_{\blacktriangleright}})\,d_{2}(h_{2_{\blacktriangleright}1_{\blacktriangleright}})\,d_{3}(h_{2_{\blacktriangleright}2_{\blacktriangleright}})\\
&= d_1(h_{1_{\blacktriangleright} 1} h_{1_{\blacktriangleright} 2 1_{\blacktriangleright}}) d_2(h_{1_{\blacktriangleright} 2 2_{\blacktriangleright}}) d_3(h_{2_{\blacktriangleright}}) \tag*{by \ref{Item:CPH3}}\\
&= d_1(h_{1_{\blacktriangleright} 1}) \underline{\varepsilon( h_{1_{\blacktriangleright} 2 1_{\blacktriangleright}}) d_2(h_{1_{\blacktriangleright} 2 2_{\blacktriangleright}})} d_3(h_{2_{\blacktriangleright}}) \\
& \quad + \uuline{\varepsilon(h_{1_{\blacktriangleright} 1}) d_1(h_{1_{\blacktriangleright} 2 1_{\blacktriangleright}}) d_2(h_{1_{\blacktriangleright} 2 2_{\blacktriangleright}})} d_3(h_{2_{\blacktriangleright}}) \tag*{since $d_1$ is k-derivation}\\
&= d_1(h_{1_{\blacktriangleright} 1}) d_2(h_{1_{\blacktriangleright} 2}) d_3(h_{2_{\blacktriangleright}})\\
& \quad + d_1(h_{1_{\blacktriangleright} 1_{\blacktriangleright}}) d_2(h_{1_{\blacktriangleright} 2_{\blacktriangleright}}) d_3(h_{2_{\blacktriangleright}}). \tag*{by Eqs.\eqref{CPH prop 1'} and \eqref{CPH counit 1}}
\end{align*}
similarly, we have
\[
\begin{aligned}
&d_{2}(h_{1_{\blacktriangleright}})\,d_{1}(h_{2_{\blacktriangleright}1_{\blacktriangleright}})\,d_{3}(h_{2_{\blacktriangleright}2_{\blacktriangleright}}) \\
&= d_{2}(h_{1_{\blacktriangleright}1})\,d_{1}(h_{1_{\blacktriangleright}2})\,d_{3}(h_{2_{\blacktriangleright}})
+ d_{2}(h_{1_{\blacktriangleright}1_{\blacktriangleright}})\,d_{1}(h_{1_{\blacktriangleright}2_{\blacktriangleright}})\,d_{3}(h_{2_{\blacktriangleright}}).
\end{aligned}
\]
Hence $\text{RHS}(h) =\text{LHS}(h)$.
Since \(h\in H\) was arbitrary, the desired identity holds.
\end{proof}

\begin{theorem}
\label{Thm the post-Lie alg of affine post-group}
Let $(G,m,\rhd)$ be an affine post-group over a field $k$.
Then its Lie algebra $\mathfrak g = \operatorname{Der}_{k, \varepsilon}(H, k)$,
equipped with the Lie bracket $[\cdot, \cdot]$ and the operation
$\triangleright$ defined above, forms a post-Lie algebra:
\[
d_1 \triangleright [d_2, d_3]
=
[d_1 \triangleright d_2, d_3]
+
[d_2, d_1 \triangleright d_3],
\]
\[
([d_1, d_2] + d_1 \triangleright d_2 - d_2 \triangleright d_1) \triangleright d_3
=
d_1 \triangleright (d_2 \triangleright d_3)
-
d_2 \triangleright (d_1 \triangleright d_3),
\]
where $d_1, d_2, d_3 \in \mathfrak{g}$.
\end{theorem}

\begin{definition}
    The post-Lie algebra $(\mathfrak{g}, [\cdot, \cdot], \triangleright)$ constructed in Theorem \ref{Thm the post-Lie alg of affine post-group} is called the {\bf post-Lie algebra of $(G,m,\rhd)$} and denoted by \(\mathcal{PL}(G,m,\rhd)\).
\end{definition}


\begin{theorem}
\label{Thm functor affine post group to postlie}
The assignment \((G,m,\rhd) \mapsto \mathcal{PL}(G,m,\rhd)\) defines a functor
\[
\mathcal{PL}\colon \mathsf{AffPGrp}_k \to \mathsf{PLie}_k
\]
from the category of affine post-group schemes to the category of post-Lie algebras.
\end{theorem}
\begin{proof}
Let $\mathcal{PL}(G,m,\rhd)=(\mathfrak{g}, [\cdot, \cdot], \triangleright)$ and $\mathcal{PL}(G',m',\rhd')=(\mathfrak{g}', [\cdot, \cdot]', \triangleright')$.
It remains to verify the functoriality on morphisms. Let \(\Phi\colon  (G, m, \rhd) \to (G', m', \rhd')\) be a morphism of affine post-groups and $\varphi\colon  (H', \blacktriangleright') \to (H, \blacktriangleright)$ be a morphism of co-post-Hopf algebras on the coordinate rings by Theorem \ref{Thm AffPGrp simeq CoHopf opp}. Then we define, for each \(d \in \mathfrak{g}\),
\[
\mathcal{PL}(\Phi)(d) \coloneq  d \circ \varphi,
\]
which is the differential of \(\Phi\). Therefore, $\mathcal{PL}(\Phi)$ is a Lie algebra homomorphism from $(\mathfrak{g},[\cdot, \cdot])$ to $(\mathfrak{g', [\cdot, \cdot]'})$. It remains to check that $\mathcal{PL}(\Phi)$ preserves the post-Lie multiplication. By the definition of co-post-Hopf algebra homomorphisms, we obtain
\begin{align*}
\mathcal{PL}(\Phi) (d_1 \triangleright d_2)&= (d_1 \triangleright d_2) \circ \varphi \\
&= \mu_k \circ (d_1 \otimes d_2) \circ \blacktriangleright \circ \varphi\\
&=\mu_k \circ (d_1 \otimes d_2)\circ (\varphi \otimes \varphi) \circ \blacktriangleright' \\
&=\mu_k \circ \big((d_1 \circ \varphi) \otimes (d_2 \circ \varphi) \big) \circ \blacktriangleright' \\
&=(d_1 \circ \varphi) \triangleright' (d_2 \circ \varphi) \\
&= \mathcal{PL}(\Phi)(d_1) \triangleright' \mathcal{PL}(\Phi)(d_2),
\end{align*}
which implies $\mathcal{PL}(\Phi)$ is a homomorphism of post-Lie algebras. Thus, $\mathcal{PL}(-)$ is a covariant functor.
\end{proof}

\begin{remark}
    When $(G,m ,\rhd)$ is an affine pre-group, $\mathcal{PL}(G,m,\rhd)$ is a commutative post-Lie algebra, i.e., a pre-Lie algebra. The converse is not true in general; see \cite[Section~10.i]{Milne17} for an example of Chevalley.   
    Combining Theorem \ref{Thm AffPGrp=AffSBr}, we obtain a functor $\mathcal{PL} \circ \mathcal{P}\colon \mathsf{AffSBr}_k \to \mathsf{PLie}_k$ which sends the affine (skew) braces introduced by Angiono, Galindo and Vendramin \cite{AGV17} to post-Lie algebras.
\end{remark}

Let $(G,m^{\rhd})$ be the sub-adjacent affine group of $(G,m,\rhd)$, whose coordinate ring $\mathcal{O}(G,m^{\rhd})= (H, \mu, \eta, \Delta_\blacktriangleright, \varepsilon, S_{\blacktriangleright})$ is the sub-adjacent Hopf algebra of $\mathcal{O}(G,m,\rhd)$. 
The underlying vector spaces of $\Lie(G,m,\rhd)$ and $\Lie(G,m^{\rhd})$ coincide, both being $\operatorname{Der}_{k, \varepsilon}(H,k)$.
Therefore, we can naturally define $\Lie(G,m^{\rhd})$ as $(\mathfrak{g}, [\cdot, \cdot]_\blacktriangleright)$:
\begin{equation}
\label{Lie algebra of subjacent affine group def}
    [d_1, d_2]_\blacktriangleright= \mu_k \circ(d_1 \otimes d_2 - d_2 \otimes d_1) \circ \Delta_\blacktriangleright.
\end{equation}

\begin{proposition}
\label{Prop subjacent Lie alg and subjacent affine group}
Let \(\mathcal{PL}(G,m,\rhd) = (\mathfrak{g}, [\cdot, \cdot], \triangleright)\) and \(\mathfrak{g}_{\triangleright}\) be sub-adjacent Lie algebra of $(\mathfrak{g}, [\cdot, \cdot], \triangleright)$. Then, the Lie algebra of the sub-adjacent affine group \((G, m^{\rhd})\) is precisely \(\mathfrak{g}_{\triangleright}\):
\[
\Lie(G, m^{\rhd}) = \mathfrak{g}_{\triangleright}.
\]
\end{proposition}
\begin{proof}
It suffices to show that the Lie brackets on the two sides coincide, namely,
\[
[d_1, d_2]_\blacktriangleright=[d_1, d_2]_\triangleright= [d_1, d_2] + d_1 \triangleright d_2- d_2 \triangleright d_1.
\]
Let $h \in H$. Then, we have
\begin{align*}
    & \quad [d_1, d_2]_\blacktriangleright (h)\\
    &=\dashuline{d_1(h_1 h_{2 1_\blacktriangleright})} d_2(h_{2 2_\blacktriangleright}) - \dashuline{d_2(h_1 h_{2 1_\blacktriangleright})} d_1(h_{2 2_\blacktriangleright}) \tag*{by Eqs. \eqref{subjacent Hopf algebra Delta} and \eqref{Lie algebra of subjacent affine group def}}\\
    &=d_1(h_1) \underline{\varepsilon(h_{2 1_\blacktriangleright}) d_2(h_{2 2_\blacktriangleright})} + \varepsilon(h_1) d_1(h_{2 1_\blacktriangleright}) d_2(h_{2 2_\blacktriangleright}) \\
    & \quad - d_2(h_1) \underline{\varepsilon(h_{2 1_\blacktriangleright}) d_1(h_{2 2_\blacktriangleright})} - \varepsilon(h_1) d_2(h_{2 1_\blacktriangleright}) d_1(h_{2 2_\blacktriangleright}) \\
    & = \underline{d_1(h_1) d_2(h_2) - d_2(h_1) d_1(h_2)} + \uuline{\varepsilon(h_1) d_1(h_{2 1_\blacktriangleright}) d_2(h_{2 2_\blacktriangleright})} - \uuline{\varepsilon(h_1) d_2(h_{2 1_\blacktriangleright}) d_1(h_{2 2_\blacktriangleright})} \tag*{by Eq. \eqref{CPH counit 1}}\\
    &= [d_1, d_2](h) + d_1(h_{1_\blacktriangleright}) d_2(h_{2_\blacktriangleright}) - d_2(h_{1_\blacktriangleright}) d_1(h_{2_\blacktriangleright}) \tag*{by Eqs. \eqref{CPH prop 1} and \eqref{Lie bracket of k-derivation}}\\
    &= [d_1, d_2](h) + d_1 \triangleright d_2 (h) - d_2 \triangleright d_1 (h).
\end{align*}
This completes the proof.
\end{proof}

In summary, by Remark \ref{Remark sub=sub}, Theorem \ref{Thm the post-Lie alg of affine post-group} and Proposition \ref{Prop subjacent Lie alg and subjacent affine group}, we have the following commutative diagram:
\begin{equation*}
\begin{tikzcd}
\substack{\text{co-post-Hopf algebra}\\ (H,\Delta, \blacktriangleright)} \arrow[r, leftrightarrow, "\simeq"] \arrow[d, "\text{sub-adjacent}"]
&\substack{\text{affine post-group}\\ (G,m,\rhd)}
  \arrow[r, "\mathcal{PL} (-)"]
  \arrow[d, "\text{sub-adjacent}"']
&
\substack{\text{post-Lie algebra}\\ (\mathfrak{g}, [\cdot, \cdot], \triangleright) }
  \arrow[d, "\text{sub-adjacent}"]
\\
\substack{\text{commutative}\\\text{Hopf algebra}\\(H,\Delta_\blacktriangleright)} 
\arrow[r, leftrightarrow, "\simeq"]
&\substack{\text{affine group}\\ (G,m^{\rhd})}
  \arrow[r, "\Lie(-)"]
&
\substack{\text{Lie algebra} \\ (\mathfrak{g},[\cdot, \cdot]_\triangleright)}.
\end{tikzcd}
\end{equation*}

The following example shows a non-trivial post-Lie algebra structure that cannot be detected from $k$-points of $\alpha_p \times \alpha_p$.
\begin{example}
\label{ex non-reduced example}
Let $k$ be a field of characteristic $p>0$, let $\alpha_p$ be the Frobenius kernel of $\mathbb G_a$, and consider the affine group scheme $G\coloneq\alpha_p \times \alpha_p$, whose coordinate ring is $\mathcal{O}(G)=k[X,Y]/(X^p,Y^p)$. Namely, for every $R\in\mathsf{Alg}_k$, $(G(R), \boxplus_R)$ is a group, where
\[
  G(R)=\{(a,b)\in R^{\times 2}:a^p=b^p=0\} \quad \text{and} \quad (a,b) \boxplus_R (c,d)=(a+c, b+d).
\]
Define $\rhd_R\colon  G(R) \times G(R) \to G(R)$ by
\[
  (a,b)\rhd_R(c,d)=(c,d+ac).
\]
This is well-defined since $(d+ac)^p=d^p+a^pc^p=0$. The maps $\rhd_R$ are natural in $R$, since algebra homomorphisms preserve addition and multiplication.
For $u=(a,b)\in G(R)$, the map
\[
  L^{\rhd_R}_u(c,d)=(c,d+ac)
\]
is an automorphism of $(G(R),\boxplus_R)$, with $(L^{\rhd_R}_u)^{-1}(c,d)=(c,d-ac)$.
For $v=(c,d)$ and $w=(r,s)$ in $G(R)$,
\[
  \bigl(u \boxplus_R(u\rhd_R v)\bigr)\rhd_R w
  =(r,s+ (a+c)r)
  =u\rhd_R(v\rhd_R w).
\]
Thus $(G,\boxplus,\rhd)$ is an affine post-group scheme. Since $(G, \boxplus)$ is commutative, it is an affine pre-group scheme. By Proposition \ref{Prop sub-adjacent affine group of affine postgroup}, its sub-adjacent affine group is $(G, \boxplus^{\rhd})$, where $(a,b)\boxplus^{\rhd}_R(c,d)=(a+c,b+d+ac)$.
Here, the coordinate ring $\mathcal{O}(G,\boxplus, \rhd)$ is a co-post-Hopf algebra with structure maps
\[
\begin{gathered}
  \Delta(X) = X \otimes 1 + 1 \otimes X,\qquad \Delta(Y) = Y \otimes 1 + 1 \otimes Y,\\
  \varepsilon(X) = \varepsilon(Y) = 0,\qquad S(X) = -X,\qquad S(Y) = -Y,\\
  \blacktriangleright(X) = 1 \otimes X,\qquad \blacktriangleright(Y) = 1 \otimes Y+X \otimes X,
\end{gathered}
\]
where $X,Y$ also denote their residue classes.
The Lie algebra $\mathfrak g=\operatorname{Der}_{k,\varepsilon}(\mathcal{O}(G),k)$ has basis $d_1,d_2$, which extract the coefficients of $X,Y$, respectively. Thus
\[
  d_1(X)=1,\quad d_1(Y)=0, \qquad d_2(X)=0,\quad d_2(Y)=1.
\]
Its Lie bracket $[-,-]_{\g}$ is zero, and $d_i\triangleright d_j= \mu_k \circ (d_i\otimes d_j) \circ \blacktriangleright$ gives
\[
  d_1\triangleright d_1=d_2, \qquad d_1\triangleright d_2=d_2\triangleright d_1=d_2\triangleright d_2=0.
\]
The sub-adjacent Lie bracket $[-,-]_{\g_\triangleright}$ also vanishes, since $(G, \boxplus^{\rhd})$ is abelian.
\end{example}

\begin{remark}
The scheme $G= \alpha_p \times \alpha_p$ is non-reduced and $G(k)=\{(0,0)\}$.
So the non-zero post-Lie product $\triangleright$ cannot be detected from $(G(k),\boxplus_k, \rhd_k)=\{ (0,0)\}$.
\end{remark}

\subsection{On affine Rota--Baxter group schemes}
In this subsection, we recall affine Rota--Baxter group schemes and the properties needed below (Subsection \ref{Sec Rota-Baxter Lie of aff RB group }), with the main goal of preparing for the differentiation problem.

\subsubsection{Rota--Baxter groups}
First, we recall the definition and examples of Rota--Baxter groups, following \cite{GLS21}.

\begin{definition}
A \textbf{Rota--Baxter group of weight \(1\)}, or simply a \textbf{Rota--Baxter group}, is a pair \((G,\mathfrak{B})\) consisting of a group \(G\) and a map \(\mathfrak{B}:G\to G\) satisfying
\[
 \mathfrak{B}(g)\mathfrak{B}(h)
 =
 \mathfrak{B}\bigl(g\,\mathfrak{B}(g)\,h\,\mathfrak{B}(g)^{-1}\bigr)
 \qquad \forall g,h\in G.
\]
The map $\mathfrak{B}$ is called the Rota--Baxter operator on $G$.
\end{definition}

\begin{definition}
    A \textbf{homomorphism} of Rota--Baxter groups $\varphi\colon  (G,\mathfrak{B}) \to (G', \mathfrak{B}')$ is a group homomorphism such that $\varphi \circ \mathfrak{B}= \mathfrak{B}' \circ \varphi$. Denote by $\mathsf{RBGrp}$ the category of Rota--Baxter groups.
\end{definition}

\begin{example}
Let $G$ be a group.
\begin{enumerate}
  \item The inverse map $(-)^{-1}:g \mapsto g^{-1}$ is a Rota--Baxter operator on $G$. 
  \item Given an exact factorization $G=HL$, a map $\mathfrak{B}:G \to G$ defined by $\mathfrak{B}(hl)=l^{-1}$ is a Rota--Baxter operator on $G$.
\end{enumerate}
\end{example}

\subsubsection{Affine Rota--Baxter group schemes}
 
\begin{definition}{\rm (\cite{MZL26})}
\label{def affine RB group}
    An \textbf{affine Rota--Baxter group scheme over $k$} (or simply an \textbf{affine Rota--Baxter group}) is a representable functor $G\colon \mathsf{Alg}_k \to \mathsf{Set}$ together with two natural transformations $m\colon G\times G \to G$ and $B:G \to G$ such that, for all $R$ in $\mathsf{Alg}_k$, 
    \[
        m_R\colon G(R)\times G(R) \to G(R) \quad \text{and} \quad B_R\colon G(R) \to G(R)
    \]
   make $\big(G(R),m_R,B_R \big)$ into a Rota--Baxter group. If $G$ is represented by a finitely generated $k$-algebra, then it is called an \textbf{affine algebraic Rota--Baxter group}.
\end{definition}

\begin{remark}
An affine Rota--Baxter group scheme is a Rota--Baxter group
object in the category of affine schemes. Equivalently, it is a
functor
\[
  G\colon \mathsf{Alg}_k \to \mathsf{RBGrp}
\]
whose underlying set-valued functor is representable. For more details, see \cite[Theorem 3.9]{MZL26}.
\end{remark}

\begin{definition}
    A \textbf{homomorphism} of affine Rota--Baxter groups \(\Phi\colon  (G, B) \to (G', B')\) is a natural transformation $\Phi\colon  G \to G'$ such that for every $R \in \mathsf{Alg}_k$, the map
\[
\Phi_R\colon  \big( G(R),B_R \big) \to \big( G'(R),B'_R \big)
\]
is a homomorphism of Rota--Baxter groups.
Denote by $\mathsf{AffRBGrp}_k$ the category of affine Rota--Baxter groups over $k$.
\end{definition}

\begin{example}{\rm (\cite{MZL26})}
Consider the general linear group scheme ${\rm GL}_n$. Define \(B_R\colon  {\rm GL}_n(R) \to {\rm GL}_n(R)\) by $B_R(M)=M^{-1}$. One can check that $B:{\rm GL}_n \to {\rm GL}_n$ is a natural transformation. Then \((\mathrm{GL}_n, B)\) is an affine algebraic Rota--Baxter group.
\end{example}

\begin{proposition}\cite{MZL26}\label{prop:bee}
    Let $(G,m,e,\operatorname{inv},B)$ be an affine Rota--Baxter group, then we have
    \begin{equation}
    \label{affine RB group Be=e}
        B \circ e= e.
    \end{equation}
\end{proposition}

\begin{proposition}{\rm \cite{MZL26}}
Any affine Rota--Baxter group $(G,m,e,\operatorname{inv},B)$ admits a new affine group structure. Define $m^B$ and $\operatorname{inv}^B$ by the following identities respectively:
\[
m^{B}\coloneq m^{(3)} \circ \big(\operatorname{id} \times B \times \operatorname{id} \times (\operatorname{inv} \circ B)\big) \circ  \langle \pi_1,\pi_1,\pi_2,\pi_1\rangle,
\]
\[
\operatorname{inv}^{B}\coloneq m^{(2)} \circ \big((\operatorname{inv} \circ B) \times \operatorname{inv} \times B \big) \circ \langle \operatorname{id}, \operatorname{id}, \operatorname{id} \rangle.
\] 
Then $(G,m^B,e,\operatorname{inv}^B)$ is an affine group, called the \textbf{descendant affine group}. Moreover, \((G, m, m^B)\) is an affine skew brace. 
\end{proposition}

Applying Proposition \ref{Prop affine skew brace give affine postgroup} on the affine skew brace \((G, m, m^B)\), we obtain
\begin{corollary}\label{cor:gmmb}
Let $(G,m,e,\operatorname{inv},B)$ be an affine Rota--Baxter group and \((G, m, m^B)\) be as above. Define a morphism $\rhd:G \times G \to G$ by the commutative diagram:
\begin{equation*}
    \begin{tikzcd}[column sep=2cm]
       G^{\times 2} \ar[r, "\rhd"] \ar[d,"{\langle \operatorname{id},\operatorname{id}\rangle \times \operatorname{id}}"] & G \\
       G^{\times 3} \ar[r, "{\operatorname{inv} \times m^{B}}"] & G^{\times 2} \ar[u, "{m}"].
    \end{tikzcd}
\end{equation*}
Then $(G, m, \rhd)$ is an affine post-group.
\end{corollary}

\subsubsection{Co-Rota--Baxter Hopf algebras}
Co-Rota--Baxter Hopf algebras are precisely the coordinate rings of affine Rota--Baxter groups. Let \((H, \mu, \eta, \Delta, \varepsilon, S)\) be a commutative Hopf algebra.  

\begin{definition}[\cite{Zheng24}]
An algebra map \(\mathcal{B} : H \to H\) is called a \textbf{Rota--Baxter co-operator} on \(H\) if for all \(x \in H\),
\begin{equation}
\label{CRBH definition}
    \mathcal{B}(x_1) \otimes \mathcal{B}(x_2) = \mathcal{B}(x)_1\mathcal{B}(\mathcal{B}(x)_2S(\mathcal{B}(x)_4)) \otimes \mathcal{B}(x)_3.
\end{equation}
Here we use Sweedler notation: $\Delta(x)=x_1 \otimes x_2$ and $\Delta^{(3)}(x)=x_1 \otimes x_2 \otimes x_3 \otimes x_4$. By a \textbf{co-Rota--Baxter Hopf algebra} we mean a pair \((H, \mathcal{B})\) consisting of a commutative Hopf algebra \(H\) and a Rota--Baxter co-operator \(\mathcal{B}\) on \(H\). 
\end{definition}

\begin{definition}
    A \textbf{homomorphism} \(\varphi\colon  (H, \mathcal{B}_H) \to (K, \mathcal{B}_K)\) of co-Rota--Baxter Hopf algebras is a Hopf algebra homomorphism satisfying \(\varphi \circ \mathcal{B}_H =\mathcal{B}_K \circ \varphi\).
Thus, co-RB Hopf algebras over \(k\) together with their homomorphisms form a category, denoted by \(\mathsf{CoRBHopf}_k\).
\end{definition}

\begin{example}
Every finite Rota--Baxter group $(G,B)$ gives rise to a co-Rota--Baxter Hopf algebra $(k^G, B^*)$, where $k^G\coloneq  \{f:G \to k\}$ is the algebra of functions over $G$ and
$B^*\colon k^G \to k^G, f \to f \circ B$.
\end{example}

\begin{theorem}[\cite{MZL26}]
\label{THm affine RB group=coRBHopf}
     Let $(H,\mu,\eta ,\Delta, \varepsilon, S)$ be a commutative Hopf algebra. If $(H,\mathcal{B})$ is co-Rota--Baxter Hopf algebra, then $(h^H,h^\Delta,h^{\mathcal{B}})$ is an affine Rota--Baxter group. Conversely, if $(G,m,B)$ is an affine Rota--Baxter group, $\mathcal{O}(G,m,B)$ is a co-Rota--Baxter Hopf algebra. Moreover, the correspondence is functorial and yields a contravariant equivalence of categories:
\[
\mathsf{CoRBHopf}^{\mathrm{opp}}_k \simeq \mathsf{AffRBGrp}_k.
\]
\end{theorem}
\begin{proof}
    See \cite[Theorems 3.18 and 3.19]{MZL26}.
\end{proof}

\begin{corollary}
    Let $(H,\mathcal{B})$ be a co-Rota--Baxter Hopf algebra, then we have
    \begin{equation}
    \label{CRBH eB=e}
        \varepsilon \circ \mathcal{B}=\varepsilon.
    \end{equation}
\end{corollary}
\begin{proof}
    By Theorem \ref{THm affine RB group=coRBHopf}, let $(G,m,B)$ be the affine Rota--Baxter group corresponding to $(H,\mathcal{B})$. Using the Yoneda lemma, Eq. \eqref{affine RB group Be=e} gives Eq. \eqref{CRBH eB=e}.
\end{proof}

A correspondence also exists between the descendant affine group and the descendant Hopf algebra.

\begin{remark}
Any co-Rota--Baxter Hopf algebra $(H, \mathcal{B})$ gives rise to a new commutative Hopf algebra $H_{\mathcal{B}}\coloneq (H, \mu, \eta, \Delta_{\mathcal{B}}, \varepsilon, S_{\mathcal{B}})$, called the {\bf descendant Hopf algebra} of $(H, \mathcal{B})$, where
\begin{equation}
\label{descendant comul}
    \Delta_{\mathcal{B}}(x) = x_1 \mathcal{B}(x_2) \mathcal{B}(S(x_4)) \otimes x_3,
\end{equation}
\begin{equation*}
    S_{\mathcal{B}}(x)=\mathcal{B}(S(x_1)) S(x_2) \mathcal{B}(x_3).
\end{equation*}
\end{remark}

\begin{proposition}\cite{MZL26}
\label{Prop descendant affine group and descendant Hopf algebra}
    Let $(H, \mathcal{B})$ be the coordinate ring of the affine Rota--Baxter group $(G,m,e,\operatorname{inv},B)$. Then, the descendant Hopf algebra $H_{\mathcal{B}}$ of $(H,\mathcal{B})$ is the coordinate ring of the descendant affine group $(G,m^{B},e,\operatorname{inv}^{B})$ of $(G,m,e,\operatorname{inv},B)$.
\end{proposition}

\begin{remark}
Since the isomorphism between the categories of affine skew braces and affine post-groups has already been established, the relations between affine Rota–Baxter groups and affine post-groups—and likewise between their coordinate rings—are immediate; we therefore do not elaborate on them here.
\end{remark}

\subsection{The Rota--Baxter Lie algebra of an affine Rota--Baxter group scheme}\label{Sec Rota-Baxter Lie of aff RB group }
We begin by constructing a Rota--Baxter Lie algebra (of weight 1) from an affine Rota--Baxter group scheme. The result is then generalized to Rota–Baxter Lie algebras (of arbitrary weight).

Let $(G,m,B)$ be an affine Rota--Baxter group scheme over $k$ and $\mathcal{O}(G,m,B)\coloneq (H,\mu,\eta, \Delta, \varepsilon, S, \mathcal{B})$ be its coordinate ring.
Let $\mathfrak g = \operatorname{Der}_{k,\varepsilon}(H,k)$ be the Lie algebra of the underlying affine group scheme $(G,m)$. We define a linear map
\[
\mathscr{B}\colon \mathfrak g \to \mathfrak g, \quad d \mapsto \mathscr{B}(d)\coloneq d \circ \mathcal{B}.
\]
The linearity of $\mathscr B$ is immediate from the definition. Then, the well-definedness of $\mathscr B$ is guaranteed by the following proposition.

\begin{proposition}
\label{Prop B(d) is again a k-derivation}
For any $d \in \mathfrak g$, we have $\mathscr B(d) \in \mathfrak g$.
\end{proposition}
\begin{proof}
For $x,y\in H$, we have,
\begin{align*}
\mathscr B(d)(xy) &= d \mathcal{B}(xy) \\
&= d \big(\mathcal{B}(x) \mathcal{B}(y)\big) \tag*{since $\mathcal B$ is algebra homomorphism} \\
&= d \mathcal{B}(x) \underline{\varepsilon \mathcal{B}(y)} + \underline{\varepsilon \mathcal{B}(x)} d\mathcal{B}(y) \tag*{since $d$ is a $k$-derivation} \\
&=d \mathcal{B}(x) \varepsilon(y) + \varepsilon(x) d \mathcal{B}(y) \tag*{by Eq. \eqref{CRBH eB=e}}\\
&= \mathscr{B}(d)(x) \varepsilon(y) + \varepsilon(x) \mathscr{B}(d)(y).
\end{align*}
Thus, $\mathscr B(d)$ is again a $k$-derivation.
\end{proof}

\begin{lemma}
    For all $d \in \mathfrak{g}$, we have
     \begin{equation}
    \label{dBS=-dB}
        d\mathcal{B}S=-d \mathcal{B}.
    \end{equation}
\end{lemma}

\begin{proof}
    It is well-known that $d\circ S=-d$ in every commutative Hopf algebra $A$ and every $k$-derivation $d$ of $A$. By Proposition \ref{Prop B(d) is again a k-derivation}, $\mathscr B(d)=d\circ\mathcal B$ is a $k$-derivation, so we are done. 
\end{proof}

\begin{theorem}
\label{Thm Rota--Baxter Lie algebra of affine RB group}
    Let $(G,m,B)$ be an affine Rota--Baxter group over $k$. Then, the linear map $\mathscr{B}$ is a Rota--Baxter operator of weight 1 on the Lie algebra $(\mathfrak{g}, [\cdot, \cdot])$. That is, for all $d_1, d_2 \in \mathfrak{g}$, we have
    \[
    [\mathscr{B}(d_1), \mathscr{B}(d_2)]= \mathscr{B}([\mathscr{B}(d_1), d_2] + [d_1, \mathscr{B}(d_2)] + [d_1, d_2]).
    \]
    Moreover, we call $(\g, [\cdot, \cdot], \mathscr{B})$ \textbf{the Rota--Baxter Lie algebra (of weight 1) of $(G,m,B)$}, denoted by $\mathcal{RBL}_1(G,m,B)$.
\end{theorem}

\begin{proof}
     Let $x \in H$ and $y\coloneq  \mathcal{B}(x)$. Hence, we have
     \[
     \Delta^{(3)}\big(\mathcal{B}(x) \big)=\mathcal{B}(x)_1 \otimes \mathcal{B}(x)_2 \otimes \mathcal{B}(x)_3 \otimes \mathcal{B}(x)_4 = y_1 \otimes y_2 \otimes y_3 \otimes y_4.
     \] 
 First, expanding the left-hand side yields:
    \[
     [\mathscr{B}(d_1), \mathscr{B}(d_2)] (x) 
     = d_1 \mathcal{B}(x_1) d_2 \mathcal{B}(x_2)- d_2\mathcal{B}(x_1) d_1 \mathcal{B}(x_2).
     \]
 Computing the first summand of $[\mathscr{B}(d_1), \mathscr{B}(d_2)] (x)$:
 \begin{align*}
     & \quad d_1 \mathcal{B}(x_1) d_2 \mathcal{B}(x_2)\\
     &=d_1\big(\mathcal{B}(x)_1\mathcal{B}(\mathcal{B}(x)_2) \mathcal{B} S(\mathcal{B}(x)_4)\big) d_2\big(\mathcal{B}(x)_3\big) \tag*{by Eq. \eqref{CRBH definition}}\\
     &= d_1\big( y_1 \mathcal{B}(y_2) \mathcal{B} S (y_4) \big) d_2(y_3).  
 \end{align*}
 Since $d_1$ is a $k$-derivation and $\varepsilon$ is an algebra homomorphism, we further have
 \begin{align*}
     & \quad d_1 \mathcal{B}(x_1) d_2 \mathcal{B}(x_2)\\
     &=d_1(y_1) \varepsilon \mathcal{B}(y_2) \varepsilon \mathcal{B} S (y_4)  d_2(y_3) + \varepsilon (y_1) d_1 \big(\mathcal{B}(y_2) \mathcal{B} S(y_4)\big) d_2 (y_3)\\
     &= d_1(y_1) \underline{\varepsilon \mathcal{B}(y_2)} \, \underline{\varepsilon \mathcal{B} S (y_4)}  d_2(y_3) \\
     & \quad + \varepsilon (y_1) d_1 \mathcal{B}(y_2) \underline{\varepsilon\mathcal{B} S (y_4)} d_2 (y_3) + \varepsilon (y_1 ) \underline{\varepsilon \mathcal{B} (y_2)}  d_1\mathcal{B} S (y_4) d_2 (y_3) \\
     & = \underline{d_1(y_1) \varepsilon(y_2)} \, \underline{\varepsilon(y_4) d_2(y_3)} + \underline{\varepsilon(y_1) d_1 \mathcal{B}(y_2)} \, \underline{\varepsilon(y_4) d_2(y_3)} + \uuline{\varepsilon(y_1) \varepsilon(y_2)} d_1 \mathcal{B} S (y_4) \uuline{d_2(y_3)} \tag*{by Eq. \eqref{CRBH eB=e} and $\varepsilon \circ S=\varepsilon$}\\
     &=d_1(y_1) d_2(y_2) + d_1 \mathcal{B}(y_1) d_2(y_2) + \underline{d_1 \mathcal{B}S(y_2)} d_2(y_1)\\
     &=d_1(y_1) d_2(y_2) + d_1 \mathcal{B}(y_1) d_2(y_2) - d_1 \mathcal{B}(y_2) d_2(y_1) \tag*{by Eq. \eqref{dBS=-dB}}.
 \end{align*}
 Similarly,
 \[
 d_2\mathcal{B}(x_1) d_1 \mathcal{B}(x_2)=d_2(y_1) d_1(y_2) + d_2 \mathcal{B}(y_1) d_1(y_2) - d_2 \mathcal{B}(y_2) d_1(y_1).
 \]
 Therefore, we obtain
 \begin{align*}
     [\mathscr{B}(d_1), \mathscr{B}(d_2)] (x) 
    &= d_1 \mathcal{B}(y_1) d_2(y_2) - d_1 \mathcal{B}(y_2) d_2(y_1) \\
     & \quad + d_2 \mathcal{B}(y_2) d_1(y_1) - d_2 \mathcal{B}(y_1) d_1(y_2)\\
     & \quad + d_1(y_1) d_2(y_2) - d_2(y_1) d_1(y_2)\\
     &= [\mathscr{B}(d_1), d_2](y) + [d_1, \mathscr{B}(d_2)](y)+ [d_1, d_2](y)\\
     &= \mathscr{B}\big([\mathscr{B}(d_1), d_2] + [d_1, \mathscr{B}(d_2)] + [d_1, d_2]\big) (x),
 \end{align*}
 as desired.
\end{proof}

This gives a canonical and functorial assignment \((G,m,B) \mapsto \mathcal{RBL}_1 (G,m,B)\).

\begin{theorem}
\label{Thm functor affine RB group to RBlie}
    The assignment $(G,m,B) \mapsto \mathcal{RBL}_1 (G,m,B)$ defines a  functor
    \[
    \mathcal{RBL}_1\colon \mathsf{AffRBGrp}_k \to \mathsf{RBLie}_{k,1}
    \]
    from the category of affine Rota--Baxter group schemes to the category of Rota--Baxter Lie algebras (of weight 1).
\end{theorem}
\begin{proof}
Let $\Psi\colon  (G,m,B) \to (G',m',B')$ be a homomorphism of affine Rota--Baxter groups. By Theorem \ref{THm affine RB group=coRBHopf}, this induces a co-Rota--Baxter Hopf algebra homomorphism on the coordinate rings $\psi\colon  (H', \mathcal{B}') \to (H, \mathcal{B})$. Define
\[
\mathcal{RBL}_1(\Psi)\coloneq (\dd\Psi)_e=(d\mapsto d \circ \psi).
\]
Denote $\mathcal{RBL}_1(G,m,B)$ and $\mathcal{RBL}_1(G',m',B')$ as $(\mathfrak{g}, [\cdot, \cdot], \mathscr{B})$ and $(\mathfrak{g}', [\cdot, \cdot]', \mathscr{B}')$, respectively. Since $\mathcal{RBL}_1(\Psi)$ is the differential of $\Psi$, it follows that $\mathcal{RBL}_1(\Psi)\colon \mathfrak{g} \to \mathfrak{g}'$ is a Lie algebra homomorphism. Moreover, by definition, we have
\begin{align*}
    \mathcal{RBL}_1 (\Psi) \circ \mathscr{B} (d) &= \mathcal{RBL}_1 (\Psi) ( d \circ \mathcal{B})\\
    &=d \circ \mathcal{B} \circ \psi\\
    &= d \circ \psi \circ \mathscr{B}'\\
    &= \mathcal{RBL}_1(\Psi)(d) \circ \mathscr{B}'= \mathcal{B}' \circ \mathcal{RBL}_1(\Psi)(d),
\end{align*}
which implies $\mathcal{RBL}_1 (\Psi)$ is a Rota--Baxter Lie algebra homomorphism. This completes the proof.
\end{proof}

Let $(G,m^{B})$ be the descendant affine group of $(G,m,B)$, whose coordinate ring $\mathcal{O}(G,m^{B})= (H, \mu, \eta, \Delta_\mathcal{B}, \varepsilon, S_{\mathcal{B}})=H_{\mathcal{B}}$ is the descendant Hopf algebra of $\mathcal{O}(G,m,B)$. 
The underlying vector spaces of $\Lie(G,m,B)$ and $\Lie(G,m^{B})$ coincide, both being $\operatorname{Der}_{k, \varepsilon}(H,k)$.
Therefore, we can naturally define $\Lie(G,m^{B})$ as $(\mathfrak{g}, [\cdot, \cdot]_\mathcal{B})$:
\begin{equation*}
    [d_1, d_2]_\mathcal{B}= \mu_k \circ(d_1 \otimes d_2 - d_2 \otimes d_1) \circ \Delta_\mathcal{B}.
\end{equation*}

\begin{proposition}
\label{Prop descendant Lie alg and descendant affine group}
Let \(\mathcal{RBL}_1(G,m,B) = (\mathfrak{g}, [\cdot, \cdot], \mathscr{B})\) be the Rota--Baxter Lie algebra of \((G,m,B)\), and let \(\mathfrak{g}_{\mathscr{B}}\) denote its descendant Lie algebra. Then the Lie algebra of the descendant affine group \((G, m^{B})\) is precisely \(\mathfrak{g}_{\mathscr{B}}\):
\[
\Lie(G, m^{B}) = \mathfrak{g}_{\mathscr{B}}.
\]
\end{proposition}
\begin{proof}
    It suffices to show that Lie brackets on the two sides coincide, namely,
    \begin{equation}\label{eq:B-LAH}
        [d_1, d_2]_{\mathcal{B}}=[d_1, d_2]_{\mathscr{B}}=[\mathscr{B}(d_1), d_2] + [d_1, \mathscr{B}(d_2)] + [d_1, d_2].
    \end{equation}
    First, expanding the left-hand side gives:
    \begin{align*}
        & \quad [d_1, d_2]_{\mathcal{B}}(h)\\
        &=\mu_k \circ (d_1 \otimes d_2- d_2\otimes d_1) \big(h_1 \mathcal{B}(h_2) \mathcal{B}S(h_4) \otimes h_3 \big) \tag*{by Eq. \eqref{descendant comul}}\\
        &=\mu_k \circ (d_1 \otimes d_2- d_2\otimes d_1) \big(h_1 \mathcal{B}\big(h_2 S(h_4)\big) \otimes h_3 \big) \\
        &=d_1 \big(h_1 \mathcal{B}\big(h_2 S(h_4)\big) \big) d_2(h_3)- d_2\big(h_1 \mathcal{B}\big(h_2 S(h_4)\big) \big) d_1(h_3).
    \end{align*}
    Then, we compute the first summand of $[d_1, d_2]_{\mathcal{B}}(h)$:
    \[
    d_1 \big(h_1 \mathcal{B}\big(h_2 S(h_4)\big) \big) d_2(h_3) =d_1(h_1) \varepsilon \mathcal{B}\big(h_2 S(h_4) \big) d_2(h_3) + \varepsilon(h_1) d_1 \mathcal{B}\big(h_2 S(h_4) \big) d_2(h_3).
    \]
By Proposition \ref{Prop B(d) is again a k-derivation}, $d \mathcal{B}$ is a $k$-derivation. Thus,
    \begin{align*}
        & \quad d_1 \big(h_1 \mathcal{B}\big(h_2 S(h_4)\big) \big) d_2(h_3)\\
        &=d_1(h_1) \underline{\varepsilon \mathcal{B} \big(h_2 S(h_4)\big)} d_2(h_3) + \varepsilon(h_1) \big( d_1 \mathcal{B}(h_2) \uuline{\varepsilon S(h_4)} + \varepsilon(h_2) \uwave{d_1 \mathcal{B} S(h_4)} \big) d_2(h_3)\\
        &=d_1(h_1) \varepsilon \big(h_2 S(h_4)\big) d_2(h_3) + \varepsilon(h_1) \big(d_1 \mathcal{B}(h_2) \varepsilon(h_4) - \varepsilon(h_2) d_1\mathcal{B}(h_4) \big) d_2 (h_3) \tag*{by Eqs. \eqref{CRBH eB=e}, \eqref{dBS=-dB} and $\varepsilon \circ S=\varepsilon$}\\
        &=d_1(h_1) \varepsilon(h_2) \underline{\varepsilon S(h_4)} d_2(h_3) + \varepsilon(h_1) d_1 \mathcal{B}(h_2) \varepsilon(h_4) d_2(h_3) - \varepsilon(h_1) \varepsilon(h_2) d_1 \mathcal{B}(h_4) d_2(h_3)\\
        &=\underline{d_1(h_1) \varepsilon(h_2)} \, \underline{\varepsilon (h_4) d_2(h_3)} + 
        \underline{\varepsilon(h_1) d_1 \mathcal{B}(h_2)} \, \underline{\varepsilon(h_4) d_2(h_3)} - \uuline{\varepsilon(h_1) \varepsilon(h_2)} d_1 \mathcal{B}(h_4) d_2(h_3)\\
        &=d_1(h_1) d_2(h_2) + d_1 \mathcal{B} (h_1) d_2(h_2) - \underline{\varepsilon(h_1)} d_1 \mathcal{B} (h_3) \underline{d_2 (h_2)}\\
        &=d_1(h_1) d_2(h_2) + d_1 \mathcal{B} (h_1) d_2(h_2) - d_1 \mathcal{B} (h_2) d_2(h_1)\\
        &= d_1(h_1) d_2(h_2) + [d_1 \mathcal{B}, d_2](h)\\
        &=d_1(h_1) d_2(h_2) + [
        \mathscr{B}(d_1), d_2](h).
    \end{align*}
Similarly, we have
\[
d_2\big(h_1 \mathcal{B}\big(h_2 S(h_4)\big) \big) d_1(h_3)= d_2(h_1) d_1(h_2) + [\mathscr{B}(d_2), d_1](h).
\]
Finally, we obtain
\begin{align*}
    [d_1, d_2]_{\mathcal{B}}(h)&=d_1(h_1) d_2(h_2) + [
        \mathscr{B}(d_1), d_2](h) - d_2(h_1) d_1(h_2) - [\mathscr{B}(d_2), d_1](h)\\
     &=[d_1, d_2](h) + [
        \mathscr{B}(d_1), d_2](h) + [d_1, \mathscr{B}(d_2)](h)\\
     &=[d_1, d_2]_{\mathscr{B}} (h),
\end{align*}
which completes the proof.
\end{proof}

In summary, by Theorem \ref{THm affine RB group=coRBHopf}, Proposition \ref{Prop descendant affine group and descendant Hopf algebra}, Theorem \ref{Thm Rota--Baxter Lie algebra of affine RB group} and Proposition \ref{Prop descendant Lie alg and descendant affine group}, we have the following commutative diagram:
\begin{equation*}
\begin{tikzcd}[column sep=1.5cm]
\substack{\text{co-RB Hopf algebra}\\ (H,\Delta, \mathcal{B})} \arrow[r, leftrightarrow, "\simeq"] \arrow[d, "\text{descendant}"]
&\substack{\text{affine RB group}\\ (G,m,B)}
  \arrow[r, "\mathcal{RBL}_1 (-)"]
  \arrow[d, "\text{descendant}"']
&
\substack{\text{RB Lie algebra}\\ (\mathfrak{g}, [\cdot, \cdot], \mathscr{B}) }
  \arrow[d, "\text{descendant}"]
\\
\substack{\text{commutative}\\\text{Hopf algebra}\\(H,\Delta_\mathcal{B})} 
\arrow[r, leftrightarrow, "\simeq"]
&\substack{\text{affine group}\\ (G,m^{B})}
  \arrow[r, "\Lie(-)"]
&
\substack{\text{Lie algebra} \\ (\mathfrak{g},[\cdot, \cdot]_\mathscr{B})},
\end{tikzcd}
\end{equation*}
where RB stands for Rota--Baxter.

Furthermore, motivated by Remark \ref{Remark cB is also RB operator}, we generalize the above results to the Rota--Baxter Lie algebras of arbitrary weight.
For any $\lambda \in k$, define a linear map $\mathscr{B}_\lambda$ by
\[
\mathscr{B}_\lambda\coloneq  \lambda \mathscr{B}, \quad \mathscr{B}_\lambda(d)= \lambda \mathscr{B}(d),
\]
where $d \in \mathfrak{g}$. It is easy to check that $\mathscr{B}_\lambda(d)$ is again a $k$-derivation. Thus, $\mathscr{B}_\lambda$ is a linear map from $\mathfrak{g}$ to $\mathfrak{g}$.

\begin{corollary}
    Let $(G,m,B)$ be an affine Rota--Baxter group. For any $\lambda \in k$, the linear map $\mathscr{B}_\lambda$ is a Rota--Baxter operator of weight $\lambda$ on $\mathfrak{g}$. We call $(\mathfrak{g}, [\cdot, \cdot], \mathscr{B}_\lambda)$ the {\bf Rota--Baxter Lie algebra (of weight $\lambda$) of $(G,m,B)$}, denoted by $\mathcal{RBL}_\lambda(G,m,B)$.
    Moreover, the assignment defines a functor:
    \[
    \mathcal{RBL}_{\lambda}\colon \mathsf{AffRBGrp}_k \to \mathsf{RBLie}_{k,\lambda},
    \]
where $\mathsf{RBLie}_{k,\lambda}$ is the category of Rota--Baxter Lie algebras (of weight $\lambda$).
\end{corollary}
\begin{proof}
    By Theorem \ref{Thm Rota--Baxter Lie algebra of affine RB group}, the linear map $\mathscr{B}$ is a Rota--Baxter operator of weight 1 on the Lie algebra $(\mathfrak{g}, [\cdot, \cdot])$. Moreover, $\mathscr{B}_\lambda$ is a Rota--Baxter operator of weight $\lambda$ on $(\mathfrak{g},[\cdot, \cdot])$ by Remark \ref{Remark cB is also RB operator}. Furthermore, the functoriality of $\mathcal{RBL}_\lambda$ is evident.
\end{proof}

\subsection{Open problems on ``algebraic integration''}

The preceding constructions naturally lead to the following questions.
\begin{enumerate}
\item Given a post-Lie algebra
$(\mathfrak g,[\cdot, \cdot]_{\g},\triangleright)$, does there exist an affine post-group whose post-Lie algebra is isomorphic to $(\mathfrak g,[\cdot,\cdot]_{\g},\triangleright)$?
\item Given a Rota--Baxter Lie algebra $(\mathfrak g, \mathscr{B})$, does there exist an affine Rota--Baxter group whose Rota--Baxter Lie algebra is isomorphic to $(\mathfrak g,\mathscr{B})$?
\end{enumerate}

We may also ask these questions with the additional requirement that the underlying affine group $G$ be \textbf{algebraic}. Thus, the associated post-Lie and Rota--Baxter Lie algebras are finite-dimensional.
We refer to such realizations as \textbf{algebraic integrations}.

Following \cite[Section~II.3]{MilneLAG}, a Lie algebra is said to be \textbf{algebraic} if it is (isomorphic to) the Lie algebra of an affine algebraic group.
A necessary condition for an algebraic integration is that the underlying Lie algebra is algebraic. 
However, not every finite-dimensional Lie algebra is algebraic. For example, over a field $k$ of characteristic
zero, Milne recalls a five-dimensional solvable Lie algebra $\left(\g=\bigoplus_{i=1}^{5}kx_i, [\cdot, \cdot]_{\g}\right)$ due to Bourbaki that is not algebraic \cite[II,~3.42]{MilneLAG}.
It is immediate from the definitions that $\mathscr B=-\operatorname{id}$ is a Rota--Baxter operator of weight $1$, and that $x\triangleright y=-[x,y]_{\g}$ defines a post-Lie algebra structure on $\g$. Neither structure admits an algebraic integration, since the underlying Lie algebra $\g$ is not algebraic.
The preceding obstruction comes entirely from the underlying Lie algebra. To understand the role of the additional operations, we therefore ask:
\begin{enumerate}
\setcounter{enumi}{2}
\item When does a post-Lie or Rota--Baxter Lie algebra with algebraic underlying Lie algebra admit an algebraic integration?
\end{enumerate}

\section{Algebraic versus analytic constructions}\label{sec:comp}
In the Lie group setting, differentiation sends post-Lie groups to post-Lie algebras \cite{BGST24} and Rota--Baxter Lie groups to Rota--Baxter Lie algebras \cite{GLS21}; these constructions are analytic and hence restricted to \(k=\mathbb R\) or \(\mathbb C\). In this section, we extend the following theorem to post-Lie groups and Rota--Baxter Lie groups; see Theorems \ref{thm:post-lie} and \ref{thm:rb-lie}.
\begin{theorem}[{\cite[Theorem~2.1, p.~158]{MilneLAG}}]
\label{thm:lie-grps}
There is a canonical functor \( L \) from the category of real (resp. complex) algebraic groups to real (resp. complex) Lie groups, which respects Lie algebras and takes \( \mathrm{GL}_n \) to \( \mathrm{GL}_n(\mathbb{R}) \) (resp. \( \mathrm{GL}_n(\mathbb{C}) \)) with its natural structure as a Lie group. It is faithful on connected algebraic groups (all algebraic groups in the complex case).
\end{theorem}
In other words, we show that when the ground field $k$ is $\mathbb R$ or $\mathbb C$, the algebraic constructions in this paper coincide with the constructions in \cite{BGST24} and \cite{GLS21} for algebraic post-Lie groups and algebraic Rota--Baxter Lie groups. 
We mimic the proof of the analogous result for affine algebraic racks and Lie racks in previous work of the third author \cite{Ta25}.

\subsection{Adjoint representations}
Henceforth, $(G,m,e,\operatorname{inv})$ denotes an affine algebraic group over a field $k$ with Lie algebra $\g$. Also, $R$ denotes a $k$-algebra, and we form the $R$-module $\g_R\coloneq\g\otimes R=\Der_{k,\eps}(\O(G),R)$. 
The following definition aligns with the terminology and notation used for the conjugation quandle scheme $\operatorname{Conj}(G)$ in \cite{Ta25}.

\begin{definition}
    A \textbf{pointed operation} on $G$ is a natural transformation $\Theta\colon G\times G\to G$ such that $\Theta(e,g)=g$ and $\Theta(g,e)=e$ for all $k$-algebras $R$ and $g\in G(R)$. Equivalently, the corresponding $k$-algebra homomorphism $\theta\colon\O(G)\to\O(G)\otimes_k\O(G)$ satisfies the \textbf{co-fixing} and \textbf{co-fixedness} identities from \cite[Sec.\ 5.2]{Ta25}:
    \[
    (\eps\otimes\id)\circ\theta=\id,\qquad(\id\otimes\eps)\circ\theta=\eta\circ\eps.
    \]
    Given $g\in G(R)$, define
    \[
    L^\Theta_g\colon G(R)\to G(R),\quad h\mapsto \Theta(g,h),\quad \Ad^\Theta_g\coloneq(\dd L^\Theta_g)_e\colon\g_R\to\g_R.
    \]
    Given $d_1,d_2\in\g$, define
    \[
    d_1\cdot_\Theta d_2\coloneq\mu\circ(d_1\otimes d_2)\circ\theta\colon \O(G)\to k.
    \]
\end{definition}

\begin{example}\label{ex:pointed-ops}
    \phantom{}
    \begin{enumerate}[(i)]
        \item Conjugation $\operatorname{conj}(g,h)\coloneq ghg^{-1}$ is a pointed operation.
        \item If $(G,m,\rhd)$ is an affine post-group, then $\Theta\coloneq\rhd$ is a pointed operation with $\theta=\blacktriangleright$ by \eqref{CPH counit 1} and \eqref{CPH counit 2}, and $d_1\cdot_\rhd d_2=d_1\tr d_2$ by \eqref{triangleright of k-derivation}.
        \item If $(G,m,B)$ is an affine Rota–Baxter group, then Proposition \ref{prop:bee} implies that $\Theta_B\coloneq\mathrm{conj}\circ(B\times\id)$ is a pointed operation. In fact, $\Theta_B$ is the post-group operation constructed in Corollary \ref{cor:gmmb}.
    \end{enumerate}
\end{example}

The following was shown in previous work of the third author; the proof is rather technical.

\begin{proposition}[\cite{Ta25}]\label{prop:pointed-ops}
    Let $\Theta$ be a pointed operation on $G$. Then:
    \begin{enumerate}[(i)]
        \item For all $g\in G(R)$ and $d_2\in\g_R$, we have
        \[
        \Ad^\Theta_g(d_2)=\mu_R\circ(g\otimes d_2)\circ \theta.
        \]
        In particular, the assignment
        \[
        \Ad^\Theta\colon G\to\End(\g),\qquad g\mapsto\Ad^\Theta_g
        \]
        is a morphism of schemes, and $\Ad^\Theta_e=\id$.
        \item For all $d_1\in\g$, the image of $d_1$ under the differential $(\dd\Ad^\Theta)_e$ is
        \[
        (\dd\Ad^\Theta)_e(d_1)=d_1\cdot_\Theta(-).
        \]
        In particular, $d_1\cdot_\theta d_2\in\g$ for all $d_1,d_2\in\g$.
    \end{enumerate}
\end{proposition}

\begin{proof}
    See \cite[(6.3), (8.3), Lem.\ 8.1, Prop.\ 8.2]{Ta25}.
\end{proof}

\begin{definition}\label{def:adjoint-rep}
    We call the morphism of schemes $\Ad^\Theta\colon G\to\End(\g)$ constructed in Proposition \ref{prop:pointed-ops}(i) the \textbf{adjoint representation} of $G$ with respect to $\Theta$.
\end{definition}

See Remark \ref{rmk:adjoint-rep} for the reason behind our choice of terminology and notation in Definition \ref{def:adjoint-rep}.

\begin{corollary}\label{cor:pointed-ops}
    In the notation of Example \ref{ex:pointed-ops}, the following hold for all $d_1,d_2\in\g$:
    \begin{enumerate}[(i)]
        \item $d_1\cdot_{\mathrm{conj}}d_2=[d_1,d_2]$.
        \item If $(G,m,\rhd)$ is an affine post-group, then $d_1\cdot_\rhd d_2=d_1\tr d_2$.
        \item If $(G,m,B)$ is an affine Rota--Baxter group, then $d_1\cdot_{\Theta_B}d_2=[\mathscr B(d_1),d_2]$.
    \end{enumerate}
\end{corollary}

\begin{proof}
    (i): This is \cite[Rem.\ 8.3]{Ta25}.

    (ii): This is Example \ref{ex:pointed-ops}(ii).

    (iii): Since $\theta_B=(B\otimes\id)\circ\operatorname{conj}^*$, Proposition \ref{prop:pointed-ops}(i) implies that
    \[
    \Ad^{\Theta_B}=\Ad^{\operatorname{conj}}\circ B\colon G\to\End(\g).
    \]
    Due to Proposition \ref{prop:bee} and the definition of $\mathscr B$, we can apply the chain rule and Proposition \ref{prop:pointed-ops}(ii) to obtain
    \[
    d_1\cdot_{\Theta_B}(-)=(\dd\Ad^{\operatorname{conj}})_e(\mathscr B(d_1))=\mathscr B(d_1)\cdot_{\operatorname{conj}}(-)
    \]
    for all $d_1\in\g$. Applying (i) completes the proof.
\end{proof}

\begin{remark}
    By Corollary \ref{cor:pointed-ops}(ii), taking $\Theta\coloneq\rhd$ in Proposition \ref{prop:pointed-ops}(ii) recovers Proposition \ref{prop:d1d2g}.
\end{remark}

\begin{corollary}\label{cor:splitting}
    Let $(G,m,B)$ be an affine Rota--Baxter group over $k$, and let $(G,m,\rhd)$ be the induced affine post-group from Corollary \ref{cor:gmmb}. Then the post-Lie algebra $\mathcal{PL}(G,m,\rhd)$ equals the splitting post-Lie algebra $(\g,[\cdot,\cdot]_\g,\triangleright)$ of the Rota--Baxter Lie algebra $\mathcal{RBL}_1(G,m,B)$ discussed in Remark \ref{rmk:splitting}.
\end{corollary}

\begin{proof}
    By Example \ref{ex:pointed-ops}(iii), we have $\rhd=\Theta_B$, so applying Corollary \ref{cor:pointed-ops}(ii--iii) completes the proof.
\end{proof}

\begin{proposition}\label{prop:sub-adjacent}
    Let $(G,m,\rhd)$ be an affine algebraic post-group, and let $(G,m^\rhd)$ be its  sub-adjacent affine group from Proposition \ref{Prop sub-adjacent affine group of affine postgroup}. Then $\Ad^\rhd\colon (G,m^\rhd)\to\mathrm{GL}(\g)$ is a morphism of affine groups.
\end{proposition}

\begin{proof}
    Fix $g,h\in G(R)$. Given an $R$-algebra $S$, weighted associativity in $G(S)$ yields $L^\tr_{g\circ h}=L^\tr_g\circ L^\tr_h$ in $\End(G_R)$. Since all of these maps fix $e$, the chain rule yields $\Ad^\tr_{g\circ h}=\Ad^\tr_g\circ\Ad^\tr_h$. Combined with the identity $\Ad^\tr_e=\id$, this shows that $\Ad^\tr_g$ is invertible with inverse $\Ad^\tr_{\operatorname{inv}^\tr (g)}$, so $\Ad^\tr$ lands in $\mathrm{GL}(\g)$.
\end{proof}

\subsection{Algebraic and analytic comparisons}
First, we establish some notation.
Henceforth, $k$ denotes either $\mathbb R$ or $\mathbb C$. By ``smooth,'' we mean ``$C^\infty$'' if $k=\mathbb R$ and ``holomorphic'' if $k=\mathbb C$. 
Given a smooth map of manifolds $f\colon M\to N$, denote the analytic differential of $f$ at $p\in M$ by
\[
f_{*p}\colon T_pM\to T_{f(p)}N,\qquad X\mapsto\left.\frac{\dd}{\dd t}\right|_{t=0} f(\gamma(t)),
\]
where $\gamma\colon U\to M$, with $U\subset k$ an open neighborhood of $0$, is a smooth curve such that $\gamma(0)=p$ and $\gamma'(0)=X$. We chose this notation to align with \cite{BGST24}. 

If $\Phi\colon X\to Y$ is a morphism of smooth affine algebraic schemes corresponding to $\varphi\colon\O(Y)\to\O(X)$, then $\Phi_k\colon X(k)\to Y(k)$ is a polynomial map of manifolds. For all $p\in X(k)$, it is well-known that the map
\[
\iota_p\colon T_pX(k)\xrightarrow{\sim} T_pX\coloneq\Der_{k,p}(\O(X),k),\qquad \iota_p(v)(f)\coloneq (f_k)_{*p}(v)
\]
is an isomorphism of vector spaces, and the chain rule yields
\begin{equation}\label{eq:chain-rule}
    \iota_{\Phi_k(p)}\circ(\Phi_k)_{*p}=(\dd\Phi)_p\circ\iota_p,
\end{equation}
where $(\dd\Phi)_p\colon T_pX\to T_{\Phi_k(p)}Y$ is the algebro-geometric differential $d\mapsto d\circ \varphi$.

Let $(G,e,m)$ be an affine algebraic group.
Since $\operatorname{char}(k)=0$, Cartier's theorem states $G$ is smooth. Theorem \ref{thm:lie-grps} implies that $G(k)$ is a Lie group with Lie algebra $\g\an\coloneq \Lie(G(k))$, and the map $\iota_G\coloneq\iota_e\colon \g\an\xrightarrow{\sim}\g$ defined above is a Lie algebra isomorphism. 

To translate between post-Lie groups and affine post-groups, we develop the following dictionary.

\begin{definition}
    Let $\Theta$ be a pointed operation on $G$. Following \cite{kinyon} and \cite{BGST24}, define
\[
\Lambda^\Theta\colon G\to\End(\g\an),\qquad g\mapsto((L^\Theta_g)_k)_{*e}
\]
and
\[
X\cdot^{\mathrm{an}}_\Theta Y\coloneq\left.\frac{\dd}{\dd t}\right|_{t=0}\Lambda^\Theta(\exp tX)Y
\]
for all $X,Y\in\g\an$. 
\end{definition}

\begin{theorem}\label{thm:pointed-ops}
    Let $\Theta$ be a pointed operation on $G$. For all $X,Y\in\g\an$, we have
    \[
    \iota_G(X\cdot^{\mathrm{an}}_\Theta Y)=\iota_G(X)\cdot_\Theta\iota_G(Y).
    \]
\end{theorem}

\begin{proof}
    Fix $g\in G(k)$. By \eqref{eq:chain-rule} with $\Phi\coloneq L^\Theta_g$, we have
    \[
    \iota_{G}\circ\Lambda^\Theta(g)=\Ad^\Theta_k(g)\circ\iota_G.
    \]
    That is,
    \[
    \Lambda^\Theta=c_\iota\circ\Ad^\Theta_k,\quad c_\iota\colon \g\an\xrightarrow{\sim}\g,\quad c_\iota(T)\coloneq\iota_G^{-1}\circ T\circ\iota_G.
    \]
    Since $c_\iota$ is $k$-linear, the chain rule yields
    \[
    \Lambda^\Theta_{*e}=c_\iota\circ(\Ad^\Theta_k)_{*e}.
    \]
    Now, apply Theorem \ref{thm:lie-grps} to $\End(\g)$. By applying \eqref{eq:chain-rule} with $\Phi\coloneq\Ad^\Theta$, we deduce from Proposition \ref{prop:pointed-ops}(ii) that
    \[
    (\Ad^\Theta_k)_{*e}(X)=(\dd\Ad^\Theta)_e(\iota_G(X))=\iota_G(X)\cdot_\Theta(-).
    \]
    Hence,
    \[
    X\cdot^{\mathrm{an}}_\Theta Y=\Lambda^\Theta_{*e}(X)(Y)=\iota_G^{-1}(\iota_G(X)\cdot_\Theta\iota_G(Y)),
    \]
    as desired.
\end{proof}

\begin{remark}\label{rmk:adjoint-rep}
    If $\Theta=\operatorname{conj}$ as in Example \ref{ex:pointed-ops}(i), then $\Lambda^{\operatorname{conj}}(g)=\Ad_g$ in the usual sense for all $g\in G(k)$. Hence, Theorem \ref{thm:pointed-ops} and Corollary \ref{cor:pointed-ops}(i) recover the $Q=\operatorname{Conj}(G)$ case of \cite[Thm.\ 1.1]{Ta25}.
\end{remark}

\begin{remark}
    Recall from the proof of Theorem \ref{thm:pointed-ops} that $\Lambda=c_\iota\circ\Ad^\tr_k$.
    Hence, Proposition \ref{prop:sub-adjacent} is the algebro-geometric analogue to the fact that $\Lambda\colon(G(k),\circ)\to\operatorname{GL}(\g\an)$ is a Lie group homomorphism, which one deduces by applying the chain rule to \cite[Thm.\ 2.4(ii)]{BGST24}.
\end{remark}

\subsubsection{Post-Lie groups}

Recall from \cite[(56)]{BGST24} that if $(G,\cdot,\rhd)$ is a post-Lie group, then its post-Lie algebra $(\g\an,[\cdot,\cdot],\TR\an)$ is defined by $X\TR\an Y\coloneq X\cdot^{\mathrm{an}}_\rhd Y$.

\begin{theorem}\label{thm:post-lie}
    Let $k\in\{\mathbb R,\mathbb C\}$, and let $(G,m,\rhd)$ be an affine algebraic post-group over $k$. Then:
    \begin{enumerate}[(i)]
        \item $(G(k),m_k,\rhd_k)$ is a post-Lie group, and $G_\tr(k)=(G(k),m_k^\rhd)$ (see Proposition \ref{Prop sub-adjacent affine group of affine postgroup}) is its sub-adjacent Lie group.
        \item The canonical map $\iota_G\colon(\g\an,[\cdot,\cdot],\TR\an)\to\mathcal{PL}(G,m,\rhd)$ is an isomorphism of post-Lie algebras.
        \item The assignment $G\mapsto\iota_G$ is natural in the following sense: for every morphism of affine post-groups $\Phi\colon G\to H$, we have
        \[
        \iota_H\circ(\Phi_k)_{*e}=\mathcal{PL}(\Phi)\circ\iota_G.
        \]
        \item $\iota_G$ is a Lie algebra isomorphism from $\Lie(G_\rhd(k))$ to $\Lie(G_\rhd)=\g_\tr$.
    \end{enumerate}
\end{theorem}

\begin{proof}
    (i): Since $\tr$ is a morphism of schemes, applying Theorem \ref{thm:lie-grps} to $G\times G$ shows that $\tr_k$ is smooth. Thus, the claim follows from Theorem \ref{thm:lie-grps} and Proposition \ref{Prop sub-adjacent affine group of affine postgroup}.

    (ii): Take $\Theta\coloneq\rhd$ in Theorem \ref{thm:pointed-ops}. Then the claim follows from Theorem \ref{thm:lie-grps} and Corollary \ref{cor:pointed-ops}(ii).

    (iii): This is a special case of \eqref{eq:chain-rule}.

    (iv): The definition of $\iota_G$ depends only on the underlying manifold of $G(k)$ and the underlying $k$-algebra of $\O(G)$, so the claim follows from applying Theorem \ref{thm:lie-grps} to $G_\rhd$ and using Proposition \ref{Prop subjacent Lie alg and subjacent affine group}.
\end{proof}

\subsubsection{Rota--Baxter Lie groups}

Recall from \cite{GLS21} that a \textbf{Rota--Baxter Lie group} (of weight 1) is a Lie group $G$ equipped with a smooth map $\mathfrak B\colon G\to G$ such that
\[
\mathfrak B(g)\mathfrak B(h)=\mathfrak B(g\mathfrak B(g)\mathfrak B(h)\mathfrak B(g)^{-1})
\]
for all $g,h\in G$. Its \textbf{descendant Lie group} is
\[
(G,*),\qquad g*h\coloneq g\mathfrak B(g)h\mathfrak B(g)^{-1},
\]
and $\mathfrak B_{*e}$ is a Rota--Baxter operator of weight 1 on the Lie algebra $\g\an$.

\begin{theorem}\label{thm:rb-lie}
    Let $k\in\{\mathbb R,\mathbb C\}$, and let $(G,m,B)$ be an affine algebraic Rota--Baxter group over $k$. Then:
    \begin{enumerate}[(i)]
        \item $(G(k),B_k)$ is a Rota--Baxter Lie group, and $(G(k),m^B_k)$ is its descendant Lie group.
        \item For all $\lambda\in k$, the canonical map \[\iota_G\colon(\g\an,[\cdot,\cdot],\lambda(B_k)_{*e})\to\mathcal{RBL}_\lambda(G,m,B)=(\g,[\cdot,\cdot],\mathscr B_\lambda)\] is an isomorphism of Rota--Baxter Lie algebras of weight $\lambda$.
        \item The assignment $G\mapsto\iota_G$ is natural in the following sense: for all $\lambda\in k$ and for every morphism of affine Rota--Baxter groups $\Phi\colon G\to H$, we have
        \[
        \iota_H\circ(\Phi_k)_{*e}=\mathcal{RBL}_\lambda(\Phi)\circ\iota_G.
        \]
        \item $\iota_G$ is a Lie algebra isomorphism from $\Lie(G(k),m^B_k)$ to $\Lie(G,m^B)=\g_\mathscr{B}$.
    \end{enumerate}
\end{theorem}

\begin{proof}
    (i): This follows from Theorem \ref{thm:lie-grps}.

    (ii): Take $\Phi\coloneq B$ in \eqref{eq:chain-rule}. Then the claim follows from Proposition \ref{prop:bee} and Theorem \ref{thm:lie-grps}.

    (iii): This is a special case of \eqref{eq:chain-rule}.

    (iv): Thanks to Proposition \ref{Prop descendant affine group and descendant Hopf algebra}, the claim follows after applying Theorem \ref{thm:lie-grps} to $(G,m^B)$ and appealing to Proposition \ref{Prop descendant Lie alg and descendant affine group}.
\end{proof}

$\mathbf{Acknowledgements.}$ This research is supported by NSFC (W2412041, 12371029) and the National Key Research and Development Program of China (2021YFA1002000).

\end{document}